\documentclass[11pt]{article}

\usepackage{amsfonts,amssymb,amsmath,graphicx,graphics}
\usepackage[english]{babel}
 \def\botcaption#1#2{\medskip\centerline{{\scshape #1.}\kern8pt
 {\rm #2}}\bigskip}

 \makeatletter
 \@addtoreset{equation}{section}
 \makeatother

 \makeatletter
 \@addtoreset{enunciato}{section}
 \makeatother

 \newcounter{enunciato}[section]

 \newtheorem{ittheorem}{Theorem}
 \newtheorem{itlemma}{Lemma}
 \newtheorem{itproposition}{Proposition}
 \newtheorem{itdefinition}{Definition}
 \newtheorem{itassumption}{Assumption}
 \newtheorem{itremark}{Remark}
 \newtheorem{itclaim}{Claim}
 \newtheorem{itcorollary}{Corollary}

 \newenvironment{theorem}{\addtocounter{enunciato}{1}
 \begin{ittheorem}}{\end{ittheorem}}

  \newenvironment{corollary}{\addtocounter{enunciato}{1}
 \begin{itcorollary}}{\end{itcorollary}}

 \newenvironment{lemma}{\addtocounter{enunciato}{1}
 \begin{itlemma}}{\end{itlemma}}

 \newenvironment{proposition}{\addtocounter{enunciato}{1}
 \begin{itproposition}}{\end{itproposition}}

 \newenvironment{definition}{\addtocounter{enunciato}{1}
 \begin{itdefinition}}{\end{itdefinition}}

  \newenvironment{assumption}{\addtocounter{enunciato}{1}
 \begin{itassumption}}{\end{itassumption}}

 \newenvironment{remark}{\addtocounter{enunciato}{1}
 \begin{itremark}}{\end{itremark}}

 \newenvironment{claim}{\addtocounter{enunciato}{1}
 \begin{itclaim}}{\end{itclaim}}

 \newenvironment{proof}{\noindent {\bf Proof.\,}
 }{\hspace*{\fill}$\square$\medskip}

 \newcommand{\be}[1]{\begin{equation}\label{#1}}
 \newcommand{\ee}{\end{equation}}

 \newcommand{\bl}[1]{\begin{lemma}\label{#1}}
 \newcommand{\el}{\end{lemma}}

 \newcommand{\br}[1]{\begin{remark}\label{#1}}
 \newcommand{\er}{\end{remark}}

 \newcommand{\bt}[1]{\begin{theorem}\label{#1}}
 \newcommand{\et}{\end{theorem}}

 \newcommand{\bd}[1]{\begin{definition}\label{#1}}
 \newcommand{\ed}{\end{definition}}

 \newcommand{\bcl}[1]{\begin{claim}\label{#1}}
 \newcommand{\ecl}{\end{claim}}

 \newcommand{\bp}[1]{\begin{proposition}\label{#1}}
 \newcommand{\ep}{\end{proposition}}

 \newcommand{\bc}[1]{\begin{corollary}\label{#1}}
 \newcommand{\ec}{\end{corollary}}

 \newcommand{\bpr}{\begin{proof}}
 \newcommand{\epr}{\end{proof}}

 \newcommand{\bi}{\begin{itemize}}
 \newcommand{\ei}{\end{itemize}}

 \newcommand{\ben}{\begin{enumerate}}
 \newcommand{\een}{\end{enumerate}}

 \def\botcaption#1#2{\medskip\centerline{{\scshape #1.}\kern8pt
 {\rm #2}}\bigskip}

 \def \ba {\begin{array}}
 \def \ea {\end{array}}

 \def \Z {{\mathbb Z}}
 \def \R {{\mathbb R}}
 
 \def \N {{\mathbb N}}

 \def \cO {{\mathcal O}}
 \def \cW {{\mathcal W}}
 \def \cD {{\mathcal D}}
 \def \cL {{\mathcal L}}
 \def \cI {{\mathcal I}}
 \def \cJ {{\mathcal J}}
 \def \cR {{\mathcal R}}
 \def \cA {{\mathcal A}}
 
 \def \cN {{\mathcal N}}
 
 \def \cP {{\mathcal P}}

 \def \cF {{\mathcal F}}

 \def \cT {{\mathcal T}}

 \def \gep {{\varepsilon}}

 \def \DOM {{\hbox{\footnotesize\rm DOM}}}
 \def \CONE {{\hbox{\footnotesize\rm CONE}}}

\begin{document}

\title{On the localized phase of a copolymer in an emulsion:
subcritical percolation regime}

\author{\renewcommand{\thefootnote}{\arabic{footnote}}
F.\ den Hollander
\footnotemark[1]\,\,\,\footnotemark[2]
\\
\renewcommand{\thefootnote}{\arabic{footnote}}
N.\ P\'etr\'elis
\footnotemark[2]
}

\footnotetext[1]
{Mathematical Institute, Leiden University, P.O.\ Box 9512,
2300 RA Leiden, The Netherlands}\,

\footnotetext[2]
{EURANDOM, P.O.\ Box 513, 5600 MB Eindhoven, The Netherlands}

\maketitle

\begin{abstract}
The present paper is a continuation of \cite{dHP07b}. The object of interest
is a two-dimensional model of a directed copolymer, consisting of a random
concatenation of hydrophobic and hydrophilic monomers, immersed in an emulsion,
consisting of large blocks of oil and water arranged in a percolation-type fashion.
The copolymer interacts with the emulsion through an interaction Hamiltonian that
favors matches and disfavors mismatches between the monomers and the solvents, in
such a way that the interaction with the oil is stronger than with the water.

The model has two regimes, supercritical and subcritical, depending on whether
the oil blocks percolate or not. In \cite{dHP07b} we focussed on the supercritical
regime and obtained a complete description of the phase diagram, which consists
of two phases separated by a single critical curve. In the present paper we focus
on the subcritical regime and show that the phase diagram consists of four phases
separated by three critical curves meeting in two tricritical points.

\vskip 0.5truecm
\noindent
\emph{AMS} 2000 \emph{subject classifications.} 60F10, 60K37, 82B27.
{\it Acknowledgment.} NP was supported by a postdoctoral fellowship from the
Netherlands Organization for Scientific Research (grant 613.000.438).
\end{abstract}

\newpage

%%%%%%%%%%%%%%%%%%%%%%%%%%%% SECTION 1 %%%%%%%%%%%%%%%%%%%%%%%%%%%%%%%%%%%%%%%%

\section{Introduction and main results}
\label{S1}

%%%%%%%%%%%%%%%%%%%%%%%%%%%%%%%%%%%%%%%%%%%%%%%%%%%%%%%%%%%%%%%%%%%%%%%%%%%%%
\begin{figure}[hbtp]
\vspace{0.5cm}
\begin{center}
\includegraphics[scale = 0.25]{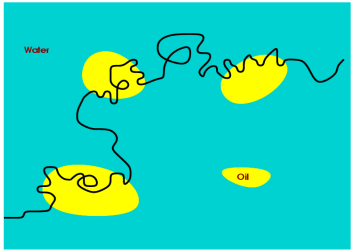}
\end{center}
\vspace{-0.5cm}
\caption{\small An undirected copolymer in an emulsion.}
\label{fig-copolemul}
\end{figure}
%%%%%%%%%%%%%%%%%%%%%%%%%%%%%%%%%%%%%%%%%%%%%%%%%%%%%%%%%%%%%%%%%%%%%%%%%%%%%

\subsection{Background}
\label{S1.1}

In the present paper we consider a two-dimensional model of a random copolymer in a random emulsion
(see Fig.~\ref{fig-copolemul}) that was introduced in den Hollander and Whittington~\cite{dHW06}.
The copolymer is a concatenation of hydrophobic and hydrophilic monomers, arranged randomly
with density $\frac12$ each. The emulsion is a collection of droplets of oil and water, arranged
randomly with density $p$, respectively, $1-p$, where $p\in(0,1)$. The configurations of the copolymer
are directed self-avoiding paths on the square lattice. The emulsion acts as a percolation-type
medium, consisting of large square blocks of oil and water, with which the copolymer interacts.
Without loss of generality we will assume that the interaction with the oil is stronger than
with the water.

In the literature most work is dedicated to a model where the solvents are separated by
a \emph{single flat infinite} interface, for which the behavior of the copolymer is the
result of an energy-entropy competition. Indeed, the copolymer prefers to match monomers
and solvents as much as possible, thereby lowering its energy, but in order to do so
it must stay close to the interface, thereby lowering its entropy. For an overview, we
refer the reader to the theses by Caravenna~\cite{C05} and P\'etr\'elis~\cite{P06}, and to
the monograph by Giacomin~\cite{G07}.

With a \emph{random interface} as considered here, the energy-entropy competition remains
relevant on the \emph{microscopic} scale of single droplets. However, it is supplemented
with the copolymer having to choose a \emph{macroscopic} strategy for the frequency at which
it visits the oil and the water droplets. For this reason, a \emph{percolation phenomenon}
arises, depending on whether the oil droplets percolate or not. Consequently, we must
distinguish between a supercritical regime $p\geq p_c$ and a subcritical regime $p< p_c$,
with $p_c$ the critical probability for directed bond percolation on the square lattice.

As was proven in den Hollander and Whittington~\cite{dHW06}, in the \emph{supercritical}
regime the copolymer undergoes a phase transition between full delocalization into the
infinite cluster of oil and partial localization near the boundary of this cluster. In
den Hollander and P\'etr\'elis \cite{dHP07b} it was shown that the critical curve separating
the two phases is strictly monotone in the interaction parameters, the phase transition is of
second order, and the free energy is infinitely differentiable off the critical curve.

The present paper is dedicated to the \emph{subcritical} regime, which turns out to be
considerably more complicated. Since the oil droplets do not percolate, even in the
delocalized phase the copolymer puts a positive fraction of its monomers in the water.
Therefore, some parts of the copolymer will lie in the water and will not localize near
the oil-water-interfaces at the same parameter values as the other parts that lie in the
oil.

\begin{figure}[htbp]
\vspace{0.5cm}
\begin{center}
\includegraphics[height=6cm]{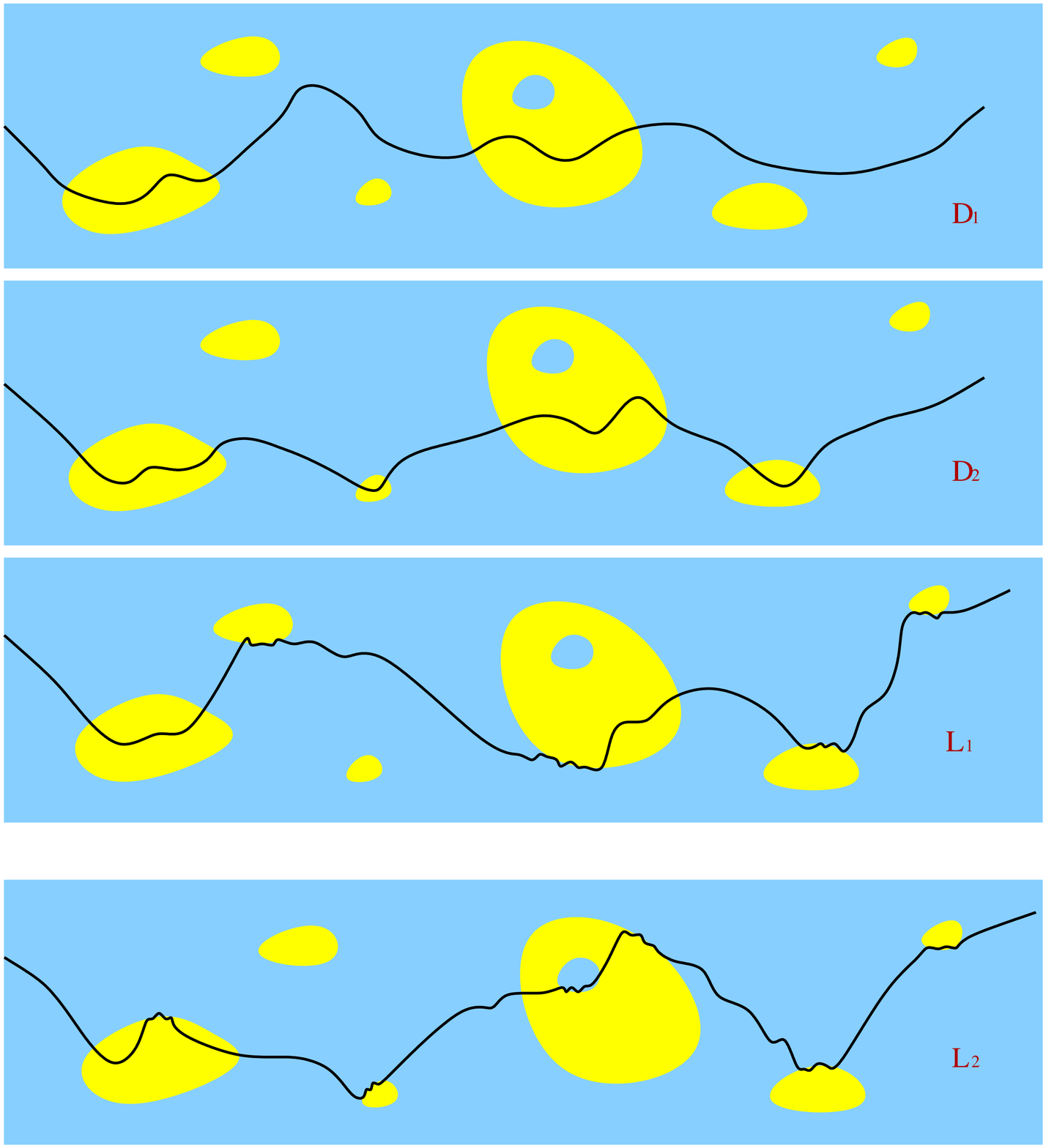}
\vspace{-0.5cm}
\end{center}
\caption{\small Typical configurations of the copolymer in each of the four phases.}
\label{fig-4phases}
\end{figure}

We show that there are four different phases (see Fig.~\ref{fig-4phases}):

\begin{itemize}
\item[(1)]
If the interaction between the two monomers and the two solvents is \emph{weak}, then the
copolymer is \emph{fully delocalized} into the oil and into the water. This means that the
copolymer crosses large clusters of oil and water alternately, without trying to follow the
interfaces between these clusters. This phase is denoted by $\cD_1$ and was investigated
in detail in \cite{dHW06}.
\item[(2)]
If the interaction strength between the \emph{hydrophobic} monomers and the two solvents is
increased, then it becomes energetically favorable for the copolymer, when it wanders around
in a water cluster, to make an excursion into the oil before returning to the water cluster.
This phase is denote by $\cD_2$ and was not noticed in \cite{dHW06}.
\item[(3)]
If, subsequently, the interaction strength between the \emph{hydrophilic} monomers and the two
solvents is increased, then it becomes energetically favorable for the copolymer, before
moving into water clusters, to follow the oil-water-interface for awhile. This phase is
denoted by $\cL_1$.
\item[(4)]
If, finally,  the interaction between the two monomers and the two solvents is \emph{strong},
then the copolymer becomes partially localized and tries to move along the oil-water
interface as much as possible. This phase is denoted by $\cL_2$.
\end{itemize}

%%%%%%%%%%%%%%%%%%%%%%%%%%%%%%%%%%%%%%%%%%%%
%%%%%%%%%%%%%%%%%%%%%%%%%%%%%%%%%%%%%%%%%%%%%

In the remainder of this section we describe the model (Section~\ref{S1.2}), recall
several key facts from \cite{dHW06} (Section~\ref{S1.3}), define and characterize
the four phases (Section~\ref{S1.4}), and prove our main results about the shape
of the critical curves and the order of the phase transitions (Section~\ref{S1.5}).

\subsection{The model}
\label{S1.2}

The \emph{randomness of the copolymer} is encoded by $\omega=(\omega_i)_{i\in \N}$,
an i.i.d.\ sequence of Bernoulli trials taking values $A$ and $B$ with probability
$\tfrac12$ each. The $i$-th monomer in the copolymer is hydrophobic when $\omega_i=A$
and hydrophilic when $\omega_i=B$. Partition $\R^2$ into square blocks of size $L_n\in\N$,
i.e.,
\be{blocks}
\R^2 = \bigcup_{x \in \Z^2} \Lambda_{L_n}(x), \qquad
\Lambda_{L_n}(x) = xL_n + (0,L_n]^2.
\ee
The \emph{randomness of the emulsion} is encoded by $\Omega=(\Omega_x)_{x\in\Z^2}$,
an i.i.d.\ field of Bernoulli trials taking values $A$ or $B$ with probability $p$,
respectively, $1-p$, where $p\in (0,1)$. The block $\Lambda_{L_n}(x)$ in the emulsion
is filled with oil when $\Omega_x=A$ and filled with water when $\Omega_x=B$.

\begin{figure}[hbtp]
\vspace{0.5cm}
\begin{center}
\includegraphics[scale = 0.4]{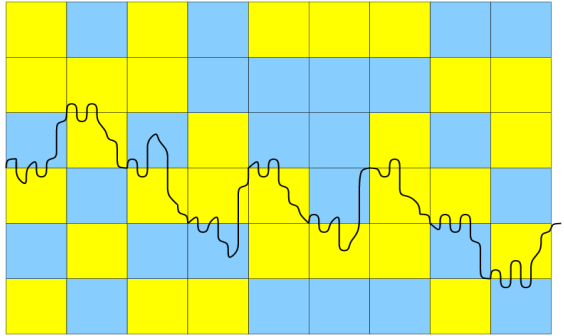}
\end{center}
\vspace{-0.5cm}
\caption{\small A directed self-avoiding path crossing blocks of oil and water
diagonally. The light-shaded blocks are oil, the dark-shaded blocks are water.
Each block is $L_n$ lattice spacings wide in both directions. The path carries
hydrophobic and hydrophilic monomers on the lattice scale, which are not indicated.}
\label{fig-copolemulblock}
\end{figure}

Let $\cW_n $ be the set of $n$-step \emph{directed self-avoiding paths} starting
at the origin and being allowed to move \emph{upwards, downwards and to the right}.
The possible configurations of the copolymer are given by a subset of $\cW_n$:
\begin{itemize}
\item
$\cW_{n,L_n} =$ the subset of $\cW_n$ consisting of those paths that enter blocks
at a corner, exit blocks at one of the two corners \emph{diagonally opposite} the
one where it entered, and in between \emph{stay confined} to the two blocks that
are seen upon entering (see Fig.~\ref{fig-copolemulblock}).
\end{itemize}
The corner restriction, which is unphysical, is put in to make the model mathematically
tractable. Despite this restriction, the model has physically relevant behavior.

%%%%%%%%%%%%%%%%%%%%%%%%%%%%%%%%%%%%%%%%%%%%%%%%%%%%%%%%%%%%%%%%%%%%%%%%%%%%%

%%%%%%%%%%%%%%%%%%%%%%%%%%%%%%%%%%%%%%%%%%%%%%%%%%%%%%%%%%%%%%%%%%%%%%%%%%%%%

Pick $\alpha,\beta \in\R$. For $\omega$, $\Omega$ and $n$ fixed, the \emph{Hamiltonian}
$H_{n,L_n}^{\omega,\Omega}(\pi)$ associated with $\pi \in \cW_{n,L_n}$ is given by
$-\alpha$ times the number of $AA$-matches plus $-\beta$ times the number of $BB$-matches.
For later convenience, we add the constant $\frac12\alpha n$, which, by the law of
large numbers for $\omega$, amounts to rewriting the Hamiltonian as
\be{Hamiltonian}
H_{n,L_n}^{\omega,\Omega}(\pi)
= \sum_{i=1}^n \Big(\alpha\, 1\left\{\omega_i=A\right\}
- \beta\, 1\left\{\omega_i=B\right\}\Big)\,
1\left\{\Omega^{L_n}_{(\pi_{i-1},\pi_i)}=B\right\},
\ee
where $(\pi_{i-1},\pi_i)$ denotes the $i$-th step in the path $\pi$ and $\Omega^{L_n}_{
(\pi_{i-1},\pi_i)}$ denotes the label of the block this step lies in. As shown in
\cite{dHW06}, Theorem 1.3.1, we may without loss of generality restrict the interaction
parameters to the cone
\be{defcone}
\CONE = \{(\alpha,\beta)\in\R^2\colon\,\alpha\geq |\beta|\}.
\ee

A path $\pi\in\cW_{n,L_n}$ can move across four different pairs of blocks. We use the
labels $k,l\in\{A,B\}$ to indicate the type of the block that is diagonally crossed,
respectively, the type of the neighboring block that is not crossed. The size $L_n$ of
the blocks in (\ref{blocks}) is assumed to satisfy the conditions
\be{Ln}
L_n \to \infty \quad \mbox{ and } \quad \frac{1}{n}\,L_n \to 0 \qquad
\text{as }  n \to \infty,
\ee
i.e., both the number of blocks visited by the copolymer and the time spent by the copolymer
in each pair of blocks tend to infinity. Consequently, the copolymer is \emph{self-averaging}
w.r.t.\ both $\Omega$ and $\omega$.

\subsection{Free energies and variational formula}
\label{S1.3}

In this section we recall several key facts about free energies from \cite{dHW06},
namely, the free energy of the copolymer near a single flat infinite interface
(Section~\ref{S1.3.1}), in a pair of neighboring blocks (Section~\ref{S1.3.2}),
respectively, in the emulsion (Section~\ref{S1.3.3}).

\subsubsection{Free energy near a single interface}
\label{S1.3.1}

Consider a copolymer in the vicinity of a \emph{single flat infinite interface}.
Suppose that the upper halfplane is oil and the lower halfplane, including the interface,
is water. For $c\geq b>0$ and $L\in \N$, let $\cW_{cL,bL}$ be the set of $cL$-step
directed self-avoiding paths from $(0,0)$ to $(bL,0)$. The entropy per step of these
paths is
\be{entpath12}
\hat{\kappa}(c/b) = \lim_{L\to \infty} \frac{1}{cL} \log |\cW_{cL,bL}|.
\ee
On this set of paths we define the Hamiltonian
\be{Zinf}
H^{\omega,\cI}_{cL}(\pi)
= \sum_{i=1}^{cL}\Big(\alpha\, 1\{\omega_i=A\} -\beta\, 1\{\omega_i=B\}\Big)\,
1\{(\pi_{i-1},\pi_i) \leq 0\},
\ee
where $(\pi_{i-1},\pi_i) \leq 0$ means that the $i$-th step lies in the lower halfplane
(as in (\ref{Hamiltonian}) we have added the constant $\frac12\alpha cL$). The associated
partition function is
\be{partfun}
Z^{\omega,\cI}_{cL,bL} = \sum_{\pi\in\cW_{cL,bL}}
\exp\left[-H^{\omega,\cI}_{cL}(\pi)\right].
\ee
It was proven in \cite{dHW06}, Lemma 2.2.1, that
\be{phidef}
\lim_{L\to \infty} \frac{1}{cL}\log Z^{\omega,\cI}_{cL,bL}
= \phi^\cI(\alpha,\beta;c/b) = \phi^{\cI}(c/b)
\qquad \omega-a.s. \mbox{ and in mean}
\ee
for some non-random function $\phi^{\cI}\colon\,[1,\infty)\to\R$.

\subsubsection{Free energy in a pair of neighboring blocks}
\label{S1.3.2}

Let $\DOM=\{(a,b)\colon\,a\geq 1+b,\,b\geq 0\}$. For $(a,b) \in \DOM$, let $\cW_{aL,bL}$
be the set of $aL$-step directed self-avoiding paths starting at $(0,0)$, ending at $(bL,L)$,
whose vertical displacement stays within $(-L,L]$ ($aL$ and $bL$ are integers). The entropy
per step of these paths is
\be{kappa}
\kappa(a,b) = \lim_{L\to\infty} \frac{1}{aL} \log |\cW_{aL,bL}|.
\ee

Explicit formulas for $\kappa$ and $\hat{\kappa}$ are given in \cite{dHW06}, Section 2.1.
These formulas are non-trivial in general, but can be used in some specific cases to perform
exact computations.

%%%%%%%%%%%%%%%%%%%%%%%%%%%%%%%%%%%%%%%%%%%%%%%%%%%%%%%%%%%%%%%%%%%%%%%%
\begin{figure}[htbp]
\vspace{0.5cm}
\begin{center}
\includegraphics[height=5cm]{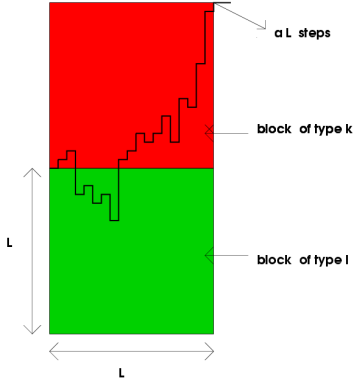}
\end{center}
\vspace{-0.5cm}
\caption{\small
Two neighbouring blocks and a piece of the path. The block that is crossed is of type $k$,
the block that appears as its neighbor is of type $l$.}
\label{fig-blockcross}
\end{figure}
%%%%%%%%%%%%%%%%%%%%%%%%%%%%%%%%%%%%%%%%%%%%%%%%%%%%%%%%%%%%%%%%%%%%%%%%%

For $k,l\in\{A,B\}$, let $\psi_{kl}$ be the \emph{quenched free energy per step} of the
directed self-avoiding path in a $kl$-block. Recall the Hamiltonian introduced in
\eqref{Hamiltonian} and for $a \geq 2$ define (see Figure \ref{fig-blockcross})
\be{fekl}
\psi_{kl}(\alpha,\beta;a) = \psi_{kl}(a)
= \lim_{L\to\infty} \frac{1}{aL} \log \sum_{\pi\in\cW_{aL,L}}
\exp\big[-H^{\omega,\Omega}_{aL,L}(\pi)\big] \qquad \omega-a.s. \mbox{ and in mean}.
\ee
As shown in \cite{dHW06}, Section 2.2, the limit exists and is non-random. For $\psi_{AA}$
and $\psi_{BB}$ explicit formulas are available, i.e.,
%(see \cite{dHW06}, Section 1.5),
\be{psiAAandBB}
\psi_{AA}(\alpha,\beta;a)=\kappa(a,1) \quad \quad \text{and} \quad \quad
\psi_{BB}(\alpha,\beta;a)=\kappa(a,1)+\frac{\beta-\alpha}{2}.
\ee
For $\psi_{AB}$
and $\psi_{BA}$  \emph{variational formulas} are available involving $\phi^{\cI}$ and $\kappa$.
To state these let, for $a\geq 2$,
\be{DOMadef}
\DOM(a) = \big\{(b,c)\in\R^2\colon\, 0\leq b\leq 1,\, c\geq b,\, a-c\geq 2-b\big\}.
\ee

\bl{psiBA} {\rm (\cite{dHW06}, Lemma 2.2.2)}\\
For all $a\geq 2$,
\be{psiinflink}
\begin{aligned}
\psi_{BA}(a) = \sup_{(b,c)\,\in\, \DOM(a)}
\frac{c\, \phi^{\cI}(c/b)+(a-c)[\frac{1}{2}(\beta-\alpha)+\kappa(a-c,1-b)]}{a}.
\end{aligned}
\ee
Moreover, $\psi_{AB}$ is given by the same expression but without the term $\frac12(\beta-\alpha)$.
\el

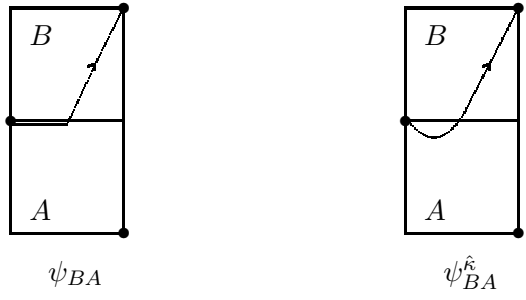
\begin{figure}[hbtp]
\vspace{0.5cm}
\begin{center}
\setlength{\unitlength}{0.25cm}
\begin{picture}(16,31)(14,-26)
\put(32,-8.5){$\psi_{BA}^{\hat{\kappa}}$}
{\thicklines
\qbezier(30,6)(33,6)(36,6)
\qbezier(36,0)(36,3)(36,6)
\qbezier(30,0)(33,0)(36,0)
\qbezier(36,-6)(36,-3)(36,0)
\qbezier(30,0)(30,3)(30,6)
\qbezier(30,0)(30,-3)(30,-6)
\qbezier(30,-6)(33,-6)(36,-6)
}
{\thicklines
\qbezier(34.15,2.9)(34.3,2.95)(34.45,3)
\qbezier(34.45,2.7)(34.45,2.85)(34.45,3)
}
\qbezier[50](33,0.2)(34.5,3.1)(36,6)
\qbezier[30](30,.2)(31.5,-2)(33,.2)
\put(31,4){$B$}
\put(31,-5.3){$A$}
\put(30,0){\circle*{.5}}
\put(36,6){\circle*{.5}}
\put(36,-6){\circle*{.5}}
\put(11,-8.5){$\psi_{BA}$}
{\thicklines
\qbezier(9,6)(12,6)(15,6)
\qbezier(15,0)(15,3)(15,6)
\qbezier(9,0)(12,0)(15,0)
\qbezier(15,-6)(15,-3)(15,0)
\qbezier(9,-0)(9,-0)(9,6)
\qbezier(9,0)(9,-3)(9,-6)
\qbezier(9,-6)(12,-6)(15,-6)
{\thicklines
\qbezier(13.15,2.9)(13.3,2.95)(13.45,3)
\qbezier(13.45,2.7)(13.45,2.85)(13.45,3)
}
}
\qbezier[50](12,-0.2)(13.5,3.1)(15,6)
\qbezier[30](9,-0.2)(10.5,-0.2)(12,-0.2)
\put(10,4){$B$}
\put(10,-5.3){$A$}
\put(9,0){\circle*{.5}}
\put(15,6){\circle*{.5}}
\put(15,-6){\circle*{.5}}
\end{picture}
\end{center}
\vspace{-4.5cm}
\caption{\small
Relevant paths for $\psi_{BA}$ and $\psi_{BA}^{\hat{\kappa}}$.}
\label{fig-subcrits2}
\end{figure}

\noindent
Similarly, we define $\psi_{BA}^{\hat{\kappa}}$ to be the free energy per step of the paths
in $\cW_{aL,L}$ that make an excursion into the $A$-block before crossing diagonally the
$B$-block, i.e.,
\be{psiinflink2}
\begin{aligned}
\psi^{\hat{\kappa}}_{BA}(a) = \sup_{(b,c)\,\in\, \DOM(a)}
\frac{c\, \hat{\kappa}(c/b)+(a-c)[\frac{1}{2}(\beta-\alpha)+\kappa(a-c,1-b)]}{a}.
\end{aligned}
\ee
Since $\hat{\kappa} \leq \phi^{\cI}$, we have $\psi_{BB}\leq\psi_{BA}^{\hat{\kappa}}
\leq\psi_{BA}$, and these inequalities are strict in some cases. The relevant paths for
(\ref{psiinflink}--\ref{psiinflink2}) are drawn in Fig.~\ref{fig-subcrits2}.

%%%%%%%%%%%%%%%%%%%%%%%%%%%%%%%%%%%%%%%%%%%%%%%%%%%%%%%%%%%%%%%%%%%%%%

%%%%%%%%%%%%%%%%%%%%%%%%%%%%%%%%%%%%%%%%%%%%%%%%%%%%%%%%%%%%%%%%%%%%%%%%%%%%%%

\br{strconcav}
{\rm
(1) As noted in \cite{dHP07b}, the strict concavity of $(a,b) \mapsto a\kappa(a,b)$ and
$\mu\mapsto\mu\hat{\kappa}(\mu)$ together with the concavity of $\mu\mapsto\mu\phi^{\cI}(\mu)$
imply that both (\ref{psiinflink}) and (\ref{psiinflink2}) have unique maximizers, which we
denote by $(\bar{b},\bar{c})$.\\
(2) In \cite{dHP07b}, we conjectured that $\mu\mapsto\mu\phi^{\cI}(\mu)$ is strictly concave.
We will need this strict concavity to prove the upper bound in Theorem \ref{th:orrder2} below.
It implies that also $a\mapsto a\psi_{BA}(a)$ and $a\mapsto a\psi_{AB}(a)$ are strictly concave.\\
(3) Since $\psi_{AA}$, $\psi_{BB}$ and $\psi_{BA}^{\hat{\kappa}}$ depend on $\alpha-\beta$ and
$a\in [2,\infty)$ only, we will sometimes write $\psi_{AA}(\alpha-\beta;a)$, $\psi_{BB}(\alpha-\beta;a)$
and $\psi_{BA}^{\hat{\kappa}}(\alpha-\beta;a)$.
}
\end{remark}

In \cite{dHW06}, Proposition 2.4.1, conditions were given under which $\bar{b},\bar{c}=0$
or $\neq 0$. Let
\be{philb4}
G(\mu,a) = \kappa(a,1) + a\partial_1\kappa(a,1) + \frac{a}{\mu}\,\partial_2\kappa(a,1)
= \frac12\left(\frac{\mu-1}{\mu}\right)\log\left(\frac{a}{a-2}\right)
+\frac{1}{\mu}\log[2(a-1)],
\ee
where $\partial_1,\partial_2$ denote the partial derivatives w.r.t.\ the first and second
argument of $\kappa(a,b)$ in (\ref{kappa}).

\bl{lemdu1}
For $a\geq 2$,
\be{mincrit}
\begin{aligned}
\psi_{AB}(a) &> \psi_{AA}(a)\quad \Longleftrightarrow \quad \sup_{\mu \geq 1}\,
\big\{\phi^{\cI}(\mu)- G(\mu,a)\big\}> 0,\\
\psi_{BA}^{\hat{\kappa}}(a)&>\psi_{BB}(a)\quad \Longleftrightarrow \quad \sup_{\mu \geq 1}\,
\big\{\hat{\kappa}(\mu)- \tfrac12(\beta-\alpha) - G(\mu,a)\big\}> 0.
\end{aligned}
\ee
\el

\subsubsection{Free energy in the emulsion}
\label{S1.3.3}

To define the \emph{quenched free energy per step} of the copolymer, we put,
for given $\omega,\Omega$ and $n$,
\be{fedef}
\begin{aligned}
f_{n,L_n}^{\omega,\Omega}
&= \frac{1}{n} \log Z_{n,L_n}^{\omega,\Omega},\\
Z_{n,L_n}^{\omega,\Omega}
&= \sum\limits_{\pi\in\cW_{n,L_n}}
\exp\left[-H_{n,L_n}^{\omega,\Omega}(\pi)\right].
\end{aligned}
\ee
As proved in \cite{dHW06}, Theorem 1.3.1,
\be{felimdef}
\lim_{n\to\infty} f_{n,L_n}^{\omega,\Omega} = f(\alpha,\beta;p)
\qquad \omega,\Omega-a.s. \mbox{ and in mean},
\ee
where, due to (\ref{Ln}), the limit is self-averaging in both $\omega$ and $\Omega$. Moreover,
$f(\alpha,\beta;p)$ can be expressed in terms of a variational formula involving the four
free energies per pair of blocks defined in (\ref{fekl}) and the frequencies at which the
copolymer visits each of these pairs of blocks on the coarse-grained block scale. To state
this variational formula, let $\mathcal{R}(p)$ be the set of $2\times 2$ matrices
$(\rho_{kl})_{k,l\in\{A,B\}}$ describing the set of possible limiting frequencies at
which $kl$-blocks are visited (see \cite{dHW06}, Section 1.3). Let $\mathcal{A}$ be the
set of $2\times 2$ matrices $(a_{kl})_{k,l\in\{A,B\}}$ such that $a_{kl}\geq 2$ for all
$k,l\in\{A,B\}$, describing the times spent by the copolymer in the $kl$-blocks on time
scale $L_n$. For $(\rho_{kl})\in \cR(p)$ and $(a_{kl})\in \cA$, we set
\be{fctV}
V\big((\rho_{kl}),(a_{kl})\big)=\frac{\sum_{kl} \rho_{kl} a_{kl} \psi_{kl}(a_{kl})}
{\sum_{kl} \rho_{kl} a_{kl}}.
\ee

\bt{feiden} {\rm (\cite{dHW06}, Theorem 1.3.1)}\\
For all $(\alpha,\beta)\in\R^2$ and $p \in (0,1)$,
\be{fevar}
f(\alpha,\beta;p) = \sup_{(\rho_{kl}) \in \cR(p)}\,\sup_{(a_{kl}) \in \cA}\,
V\big((\rho_{kl}),(a_{kl})\big).
\ee
\et

The reason why the behavior of the copolymer changes drastically at $p=p_c$ comes from the
structure of $\cR(p)$ (see Fig.~\ref{fig-rho}). For $p\geq p_c$, the set $\cR(p)$ contains
matrices $(\rho_{kl})$ satisfying $\rho_A=\rho_{AA}+\rho_{AB}=1$, i.e., the copolymer can
spend all its time inside the infinite cluster of $A$-blocks. For $p<p_c$, however, $\cR(p)$
does not contain such matrices, and this causes that the copolymer has to cross $B$-blocks
with a positive frequency. In the present paper we focus on the case $p<p_c$.

\subsection{Characterization of the four phases}
\label{S1.4}

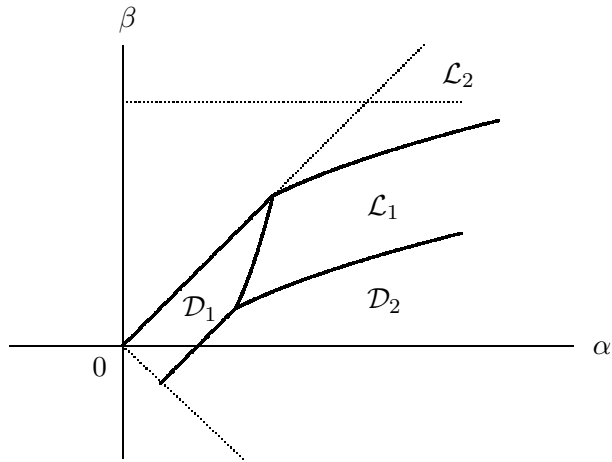
\begin{figure}[hbtp]
\begin{center}
\setlength{\unitlength}{0.5cm}
\begin{picture}(12,12)(0,-2.5)
\put(0,0){\line(12,0){12}}
\put(0,0){\line(0,8){8}}
\put(0,0){\line(0,-3){3}}
\put(0,0){\line(-3,0){3}}
{\thicklines
\qbezier(0,0)(2,2)(4,4)
\qbezier(4,4)(3.5,2)(3,1)
\qbezier(3,1)(2,0)(1,-1)
\qbezier(3,1)(5,2)(9,3)
\qbezier(4,4)(6,5)(10,6)
}
\qbezier[60](4,4)(6,6)(8,8)
\qbezier[40](0,0)(1.6,-1.5)(3.2,-3.0)
\qbezier[80](0,6.5)(4.5,6.5)(9,6.5)
\put(-.8,-.8){$0$}
\put(12.5,-0.2){$\alpha$}
\put(-0.1,8.5){$\beta$}
\put(1.6,.8){$\cD_1$}
\put(6.5,1){$\cD_2$}
\put(6.5,3.5){$\cL_1$}
\put(8.5,7){$\cL_2$}
\end{picture}
\end{center}
\caption{\small
Sketch of the phase diagram for $p<p_c$.}
\label{fig-subcrit}
\end{figure}

\begin{figure}[hbtp]
\begin{center}

\setlength{\unitlength}{0.25cm}
\begin{picture}(10,31)(12,-26)

%%% \cD_1 %%%

{\thicklines
\qbezier(0,6)(3,6)(6,6)
\qbezier(6,0)(6,3)(6,6)
\qbezier(0,0)(3,0)(6,0)
\qbezier(6,-6)(6,-3)(6,0)
\qbezier(0,0)(0,3)(0,6)
\qbezier(0,0)(0,-3)(0,-6)
\qbezier(0,-6)(3,-6)(6,-6)
}
{\thicklines
\qbezier(2.7,3)(2.85,3)(3,3)
\qbezier(3,2.7)(3,2.85)(3,3)
}
\qbezier[60](0,0)(3,3)(6,6)
\put(1,4){$A$}
\put(1,-5.3){$B$}
\put(7,-8.5){$\cD_1$}

{\thicklines
\qbezier(9,6)(12,6)(15,6)
\qbezier(15,0)(15,3)(15,6)
\qbezier(9,0)(12,0)(15,0)
\qbezier(15,-6)(15,-3)(15,0)
\qbezier(9,0)(9,3)(9,6)
\qbezier(9,0)(9,-3)(9,-6)
\qbezier(9,-6)(12,-6)(15,-6)
}
{\thicklines
\qbezier(11.7,3)(11.85,3)(12,3)
\qbezier(12,2.7)(12,2.85)(12,3)
}
\qbezier[60](9,0)(12,3)(15,6)
\put(10,4){$B$}
\put(10,-5.3){$A$}

\put(0,0){\circle*{.5}}
\put(6,6){\circle*{.5}}
\put(6,-6){\circle*{.5}}
\put(9,0){\circle*{.5}}
\put(15,6){\circle*{.5}}
\put(15,-6){\circle*{.5}}

%%% \cD_2 %%%

{\thicklines
\qbezier(21,6)(24,6)(27,6)
\qbezier(27,0)(27,3)(27,6)
\qbezier(21,0)(24,0)(27,0)
\qbezier(27,-6)(27,-3)(27,0)
\qbezier(21,0)(21,3)(21,6)
\qbezier(21,0)(21,-3)(21,-6)
\qbezier(21,-6)(24,-6)(27,-6)
}
{\thicklines
\qbezier(23.7,3)(23.85,3)(24,3)
\qbezier(24,2.7)(24,2.85)(24,3)
}
\qbezier[60](21,0)(24,3)(27,6)
\put(22,4){$A$}
\put(22,-5.3){$B$}
\put(28,-8.5){$\cD_2$}

{\thicklines
\qbezier(30,6)(33,6)(36,6)
\qbezier(36,0)(36,3)(36,6)
\qbezier(30,0)(33,0)(36,0)
\qbezier(36,-6)(36,-3)(36,0)
\qbezier(30,0)(30,3)(30,6)
\qbezier(30,0)(30,-3)(30,-6)
\qbezier(30,-6)(33,-6)(36,-6)
}
{\thicklines
\qbezier(34.15,2.9)(34.3,2.95)(34.45,3)
\qbezier(34.45,2.7)(34.45,2.85)(34.45,3)
}
\qbezier[50](33,0.2)(34.5,3.1)(36,6)
\qbezier[30](30,.2)(31.5,-2)(33,.2)
\put(31,4){$B$}
\put(31,-5.3){$A$}

\put(21,0){\circle*{.5}}
\put(27,6){\circle*{.5}}
\put(27,-6){\circle*{.5}}
\put(30,0){\circle*{.5}}
\put(36,6){\circle*{.5}}
\put(36,-6){\circle*{.5}}

%%% \cL_1 %%%

{\thicklines
\qbezier(0,-12)(3,-12)(6,-12)
\qbezier(6,-18)(6,-15)(6,-12)
\qbezier(0,-18)(3,-18)(6,-18)
\qbezier(6,-24)(6,-21)(6,-18)
\qbezier(0,-18)(0,-15)(0,-12)
\qbezier(0,-18)(0,-21)(0,-24)
\qbezier(0,-24)(3,-24)(6,-24)
}
{\thicklines
\qbezier(2.7,-15)(2.85,-15)(3,-15)
\qbezier(3,-15.3)(3,-14.85)(3,-15)
}
\qbezier[60](0,-18)(3,-15)(6,-12)
\put(1,-14){$A$}
\put(1,-23.3){$B$}
\put(7,-26.5){$\cL_1$}

{\thicklines
\qbezier(9,-12)(12,-12)(15,-12)
\qbezier(15,-18)(15,-15)(15,-12)
\qbezier(9,-18)(12,-18)(15,-18)
\qbezier(15,-24)(15,-21)(15,-18)
\qbezier(9,-18)(9,-18)(9,-12)
\qbezier(9,-18)(9,-21)(9,-24)
\qbezier(9,-24)(12,-24)(15,-24)
}
{\thicklines
\qbezier(13.15,-15.1)(13.3,-15.05)(13.45,-15)
\qbezier(13.45,-15.3)(13.45,-15.15)(13.45,-15)
}
\qbezier[50](12,-18.2)(13.5,-14.9)(15,-12)
\qbezier[30](9,-18.2)(10.5,-18.2)(12,-18.2)
\put(10,-14){$B$}
\put(10,-23.3){$A$}

\put(0,-18){\circle*{.5}}
\put(6,-12){\circle*{.5}}
\put(6,-24){\circle*{.5}}
\put(9,-18){\circle*{.5}}
\put(15,-12){\circle*{.5}}
\put(15,-24){\circle*{.5}}

%%% \cL_2 %%%

{\thicklines
\qbezier(21,-12)(24,-12)(27,-12)
\qbezier(27,-18)(27,-15)(27,-12)
\qbezier(21,-18)(24,-18)(27,-18)
\qbezier(27,-24)(27,-21)(27,-18)
\qbezier(21,-18)(21,-15)(21,-12)
\qbezier(21,-18)(21,-21)(21,-24)
\qbezier(21,-24)(24,-24)(27,-24)
}
{\thicklines
\qbezier(25.15,-15.1)(25.3,-15.05)(25.45,-15)
\qbezier(25.45,-15.3)(25.45,-15.15)(25.45,-15)
}
\qbezier[50](24,-18.2)(25.5,-14.9)(27,-12)
\qbezier[30](21,-18.2)(22.5,-18.2)(24,-18.2)
\put(22,-14){$A$}
\put(22,-23.3){$B$}
\put(28,-26.5){$\cL_2$}

{\thicklines
\qbezier(30,-12)(33,-12)(36,-12)
\qbezier(36,-18)(36,-15)(36,-12)
\qbezier(30,-18)(33,-18)(36,-18)
\qbezier(30,-18)(30,-15)(30,-12)
\qbezier(30,-18)(30,-21)(30,-24)
\qbezier(36,-18)(36,-21)(36,-24)
\qbezier(30,-24)(33,-24)(36,-24)
}
{\thicklines
\qbezier(34.15,-15.1)(34.3,-15.05)(34.45,-15)
\qbezier(34.45,-15.3)(34.45,-15.15)(34.45,-15)
}
\qbezier[50](33,-18.2)(34.5,-15.1)(36,-12)
\qbezier[30](30,-18.2)(31.5,-18.2)(33,-18.2)
\put(31,-14){$B$}
\put(31,-23.3){$A$}
\qbezier[50](33,-18.2)(34.5,-15.1)(36,-12)
\qbezier[30](30,-18.2)(31.5,-18.2)(33,-18.2)

\put(21,-18){\circle*{.5}}
\put(27,-12){\circle*{.5}}
\put(27,-24){\circle*{.5}}
\put(30,-18){\circle*{.5}}
\put(36,-12){\circle*{.5}}
\put(36,-24){\circle*{.5}}

\end{picture}

\end{center}
\caption{\small
Behavior of the copolymer inside the four block pairs containing oil and water for each
of the four phases.}
\label{fig-subcrits}
\end{figure}
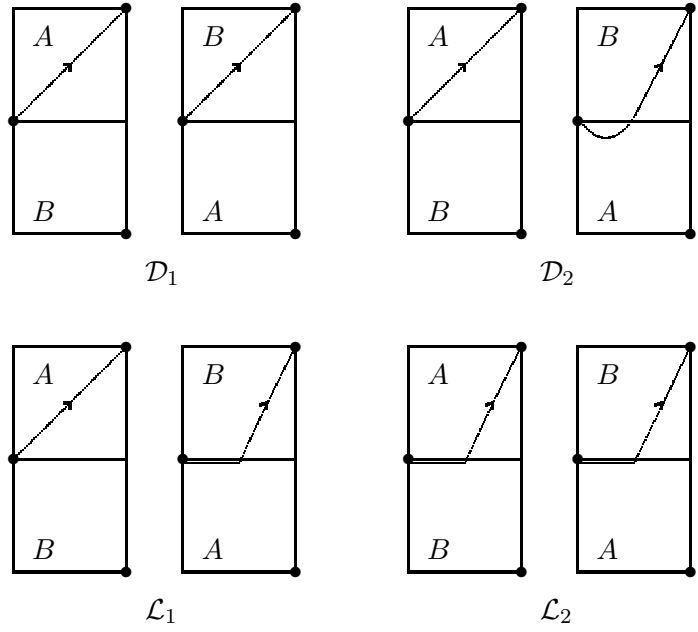

The four phases are characterized in Sections~\ref{S1.4.1}--\ref{S1.4.4}. This will involve
four free energies
\be{fefour}
f_{\cD_1} \leq f_{\cD_2} \leq f_{\cL_1} \leq f_{\cL_2} = f,
\ee
with the inequalities becoming strict successively. We will see that the phase diagram looks
like Fig.~\ref{fig-subcrit}. Furthermore, we will see that the typical path behavior in the
four phases looks like Fig.~\ref{fig-subcrits}.

%%%%%%%%%%%%%%%%%%%%%%%%%%%%%%%%%%%%%%%%%%%%%%

%%%%%%%%%%%%%%%%%%%%%%%%%%%%%%%%%%%%%%%%%%%%%%%%%%%%%%%%%%%%%%%%

%%%%%%%%%%%%%%%%%%%%%%%%%%%%%%%%%%%%%%%%%%%%%%%%%%%%%%%%%%%%%%%%%%%%%%%%

%%%%%%%%%%%%%%%%%%%%%%%%%%%%%%%%%%%%%%%%%%%%%%%%%%%%%%%%%%%%%%%%%%%%%%%%%%%

\subsubsection{The $\cD_1$-phase: $A$-delocalization and $B$-delocalization}
\label{S1.4.1}

A first region in which the free energy is analytic has been exhibited in \cite{dHW06}.
This region corresponds to full delocalization into the $A$-blocks and $B$-blocks, i.e.,
when the copolymer crosses an $AB$-block or a $BA$-block it does not spend appreciable
time near the $AB$-interface (see Fig.~\ref{fig-subcrits}). Consequently, in $\cD_1$ the
free energy depends on $\alpha-\beta$ and $p$ only, since it can be expressed in terms
of $\psi_{AA}$ and $\psi_{BB}$, which are functions of $\alpha-\beta$ (see Remark
\ref{strconcav}(3)).

\begin{figure}[hbtp]
\vspace{-1.4cm}
\begin{center}
\setlength{\unitlength}{0.45cm}
\begin{picture}(12,12)(-1,0)
\put(0,0){\line(9,0){9}}
\put(0,0){\line(0,7){7}}
{\thicklines
\qbezier(5,6)(6.25,6)(7.5,6)
\qbezier(0,0)(1,5)(5,6)
}
\qbezier[40](5,0)(5,3)(5,6)
\qbezier[40](0,6)(2.5,6)(5,6)
\qbezier[40](7.5,0)(7.5,3)(7.5,6)
\qbezier[60](0,0)(3.75,3)(7.5,6)
\put(-.6,-.8){$0$}
\put(4.8,-.8){$p_c$}
\put(5,6){\circle*{.35}}
\put(9.5,-.3){$p$}
\put(-.5,7.5){$\rho^*(p)$}
\put(-.7,5.8){$1$}
\put(7.3,-.8){$1$}
\end{picture}
\end{center}
\caption{\small
Sketch of $p\mapsto\rho^*(p)$.}
\label{fig-rho}
\end{figure}
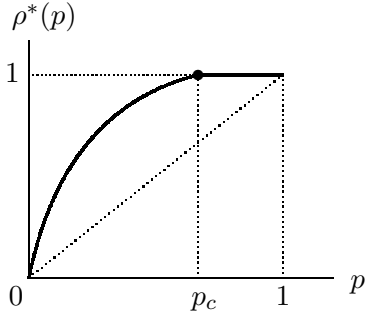

\bd{D1def}
For $p<p_c$,
\be{phaseD1}
\cD_1=\big\{(\alpha,\beta)\in \CONE\colon\,f(\alpha,\beta;p)=f_{\cD_1}(\alpha-\beta;p)\big\}
\ee
with
\be{fevarred}
f_{\cD_1}(\alpha-\beta;p) = \sup_{x\geq2,\, y\geq 2}\,
\frac{\rho^{*}(p)\, x\,\psi_{AA}(x) + [1-\rho^{*}(p)]\, y\,\psi_{BB}(y)}
{\rho^*(p) \,x + [1-\rho^*(p)]\, y},
\ee
where $\rho^*(p)$ is the maximal frequency at which the $A$-blocks can be crossed,
defined by (see Fig.~{\rm \ref{fig-rho}})
\be{rho*def}
\rho^*(p) = \max_{(\rho_{kl})\in\cR(p)} [\,\rho_{AA}+\rho_{AB}\,].
\ee
\ed

%%%%%%%%%%%%%%%%%%%%%%%%%%%%%%%%%%%%%%%%%%%%%%%%%%%%%%%%%%%%%

%%%%%%%%%%%%%%%%%%%%%%%%%%%%%%%%%%%%%%%%%%%%%%%%%%%%%%%%%%%%%%%%%

The variational formula in (\ref{fevarred}) was investigated in \cite{dHW06}, Section 2.5,
where it was found that the supremum is uniquely attained at $(\overline{x},\overline{y})$
solving the equations
\be{xysol4}
\begin{aligned}
0 &= \log 2+\rho\log(\overline{x}-2)+(1-\rho)\log(\overline{y}-2),\\
0 &= (\alpha-\beta)
+\log\left(\frac{\overline{x}(\overline{y}-2)}{\overline{y}(\overline{x}-2)}\right).
\end{aligned}
\ee
With the help of the implicit function theorem it was further proven that $f_{\cD_1}$ is
analytic on $\CONE$.

The following criteria were derived to decide whether or not $(\alpha,\beta)\in\cD_1$. The
first is a condition in terms of block pair free energies, the second in terms of the single
interface free energy.

\bp{philbex}\label{recobal} {\rm (\cite{dHW06}, Theorem 1.5.2)}
\be{D1D1cchar}
\begin{aligned}
\cD_1 &=\big\{(\alpha,\beta)\in \CONE\colon\,
\psi_{BA}(\alpha,\beta;\overline{y})=\psi_{BB}(\alpha-\beta;\overline{y})\big\},\\
\cD_1^c &=\big\{(\alpha,\beta)\in \CONE\colon\,
\psi_{BA}(\alpha,\beta;\overline{y})>\psi_{BB}(\alpha-\beta;\overline{y})\big\}.
\end{aligned}
\ee
\ep

\bc{philbextext} {\rm (\cite{dHW06}, Proposition 2.4.1 and Section 4.2.2)}
\be{D1D1ccahralt}
\begin{aligned}
\cD_1&=\Big\{(\alpha,\beta)\in \CONE\colon\, \sup_{\mu \geq 1}\,
\{\phi^{\cI}(\mu)- \tfrac12(\beta-\alpha) - G(\mu,\overline{y})\}\leq 0\Big\},\\
\cD_1^c&=\Big\{(\alpha,\beta)\in \CONE\colon\, \sup_{\mu \geq 1}\,
\{\phi^{\cI}(\mu)- \tfrac12(\beta-\alpha) - G(\mu,\overline{y})\}> 0\Big\}.
\end{aligned}
\ee
\ec

\noindent
Corollary \ref{philbextext} expresses that leaving $\cD_1$ is associated with a change
in the optimal strategy of the copolymer inside the $BA$-blocks. Namely, $(\alpha,\beta)
\in\cD_1^c$ when it is favorable for the copolymer to make an excursion into the
neighboring $A$-block before it diagonally crosses the $B$-block. This change comes
with a non-analyticity of the free energy. A first critical curve divides the phase space
into $\cD_1$ and $\cD_1^c$ (see Fig.\ \ref{fig-subcrit}).

\subsubsection{The $\cD_2$-phase: $A$-delocalization, $BA$-delocalization}
\label{S1.4.2}

Starting from $(\alpha, \beta)\in \cD_1$ with $\beta\leq 0$, we increase $\alpha$ until
it becomes energetically advantageous for the copolymer to spend some time in the $A$-solvent
when crossing a $BA$-block. It turns out that the copolymer does not localize along the
$BA$-interface, but rather crosses the interface to make a long excursion inside the $A$-block
before returning to the $B$-block to cross it diagonally (see Fig.~\ref{fig-subcrits}).

\bd{D2def}
For $p<p_c$,
\be{phaseD2}
\cD_2=\big\{(\alpha,\beta)\in\CONE\colon\,
f_{\cD_1}(\alpha-\beta;p)<f(\alpha,\beta;p)=f_{\cD_2}(\alpha-\beta;p)\big\}
\ee
with
\be{fevarredD2}
f_{\cD_2}(\alpha-\beta;p) = \sup_{x\geq 2,\, y\geq 2,\, z\geq 2}\,\sup_{\rho\in\cR(p)}
\frac{\rho_A\, x\, \psi_{AA}(x) +  \rho_{BA}\, y\,\psi^{\hat{\kappa}}_{BA}(y)
+ \rho_{BB}\, z\,\psi_{BB}(z)}{\rho_A\, x\, + \rho_{BA}\, y\,+ \rho_{BB}\, z},
\ee
where $\rho_A=\rho_{AB}+\rho_{AA}$.
\ed

\noindent
Note that $f_{\cD_2}$ depends on $\alpha-\beta$ and $p$ only, since $\psi_{AA}$, $\psi_{BB}$
and $\psi_{BA}^{\hat{\kappa}}$ are functions of $\alpha-\beta$ (see Remark \ref{strconcav}(3)).
Note also that, like (\ref{fevarred}), the variational formula in (\ref{fevarredD2}) is explicit
because we have an explicit expression for $\psi_{BA}^{\hat{\kappa}}$ via (\ref{psiinflink2})
and for $\hat{\kappa}$ and $\kappa$ via the formulas that are available from \cite{dHW06}. This
allows us to give a characterization of $\cD_2$ in terms of the block pair free energies and the
single interface free energy. For this we need a result proven in Section \ref{S2.1}, which states
that, by the strict concavity of $x\mapsto x\psi_{AA}(x)$, $y\mapsto y \psi^{\hat{\kappa}}_{BA}(y)$
and $z\mapsto z\psi_{BB}(z)$, the maximizers $(\overline{x},\overline{y},\overline{z})$ of
(\ref{fevarredD2}) are unique and do not depend on the choice of $(\rho_{kl})$ that achieves
the maximum in (\ref{fevar}).

\bp{d2char}
\be{intersec2}
\begin{aligned}
\cD_2 &= \cD_1^c\cap
\big\{(\alpha,\beta)\in \CONE \colon\;
\psi_{AB}(\overline{x})=\psi_{AA}(\overline{x})
\text{ and }
\psi_{BA}(\overline{y})=\psi_{BA}^{\hat{\kappa}}(\overline{y})\big\},\\
\cD_2^c &=\cD_1\cup
\big\{(\alpha,\beta)\in \CONE \colon\;
\psi_{AB}(\overline{x})>\psi_{AA}(\overline{x})
\text{ or }
\psi_{BA}(\overline{y})>\psi_{BA}^{\hat{\kappa}}(\overline{y})\big\}.
\end{aligned}
\ee
\ep

\bc{philbextext3}
\be{D2D2cchar}
\begin{aligned}
\cD_2 &=\cD_1^c\cap \Big\{(\alpha,\beta)\in \CONE \colon\, \sup_{\mu \geq 1}\,
\big\{\phi^{\cI}(\mu)-G(\mu,\overline{x})\big\}\leq 0 \text{ and }
\phi^{\cI}(\bar{c}/\bar{b})=\hat{\kappa}(\bar{c}/\bar{b})\Big\},\\
\cD_2^c&=\cD_1\cup\Big\{(\alpha,\beta)\in \CONE\colon\, \sup_{\mu \geq 1}\,
\big\{\phi^{\cI}(\mu)-G(\mu,\overline{x})\big\}> 0 \text{ or }
\phi^{\cI}(\bar{c}/\bar{b})>\hat{\kappa}(\bar{c}/\bar{b})\Big\},
\end{aligned}
\ee
where $(\bar{b},\bar{c})$ are the unique maximizers of the variational formula
for $\psi^{\hat{\kappa}}_{BA}(\overline{y})$ in {\rm (\ref{psiinflink2})}.
\ec

\subsubsection{The $\cL_1$-phase: $A$-delocalization, $BA$-localization}
\label{S1.4.3}

Starting from $(\alpha,\beta)\in\cD_2$, we increase $\beta$ and enter into a third phase
denoted by $\cL_1$. This phase is characterized by a partial localization along the
interface in the $BA$-blocks. The difference with the phase $\cD_2$ is that, in $\cL_1$,
the copolymer crosses the $BA$-blocks by first sticking to the interface for awhile before
crossing diagonally the $B$-block, whereas in $\cD_2$ the copolymer wanders for awhile inside
the $A$-block before crossing diagonally the $B$-block (see Fig.~\ref{fig-subcrits}). This
difference appears in the variational formula, because the free energy in the
$BA$-block is given by $\psi_{BA}$ in $\cL_1$ instead of $\psi^{\hat{\kappa}}_{BA}$
in $\cD_2$:

\bd{L1def}
For $p<p_c$,
\be{phaseL1}
\cL_1 = \big\{(\alpha,\beta)\in\CONE\colon\,
f_{\cD_2}(\alpha-\beta;p)<f(\alpha,\beta;p)=f_{\cL_1}(\alpha,\beta;p)\big\}
\ee
with
\be{fevarredL1}
f_{\cL_1}(\alpha,\beta;p) = \sup_{x\geq 2,\, y\geq 2,\, z\geq 2}\,\sup_{(\rho_{kl})\in\cR(p)}
\frac{\rho_A\, x\, \psi_{AA}(x) +  \rho_{BA}\, y\,\psi_{BA}(y) + \rho_{BB}\, z\,\psi_{BB}(z)}
{\rho_A\, z\, + \rho_{BA}\, y\,+ \rho_{BB}\, z}.
\ee
\ed

Since the strict concavity of $x\mapsto x\psi_{BA}(x)$ has not been proven (recall
Remark~\ref{strconcav}(2)), the maximizers $(\overline{x},\overline{y},\overline{z})$ of
(\ref{fevarredL1}) are not known to be unique. However, the strict concavity of $x\mapsto
x \psi_{AA}(x)$ and $z\mapsto
z \psi_{BB}(z)$ ensure that at least $\overline{x}$ and $\overline{z}$ are unique.

\bp{L1char}
\be{intersec3}
\begin{aligned}
\cL_1 &= \cD_1^c\cap\cD_2^c\cap
\big\{(\alpha,\beta)\in\CONE\colon\, \psi_{AB}(\overline{x})=\psi_{AA}(\overline{x})\big\},\\
\cL_1^c&= \cD_1\cup\cD_2\cup
\big\{(\alpha,\beta)\in\CONE\colon\, \psi_{AB}(\overline{x})>\psi_{AA}(\overline{x})\big\}.
\end{aligned}
\ee
\ep

\bc{corola3}
\be{L1charalt}
\begin{aligned}
\cL_1 &= \cD_1^c\cap\cD_2^c\cap\Big\{(\alpha,\beta)\in\CONE\colon\,
\sup_{\mu \geq 1}\,\{\phi^{\cI}(\mu) - G(\mu,\overline{x})\}\leq 0\Big\},\\
\cL_1^c &= \cD_1\cup\cD_2\cup\Big\{(\alpha,\beta)\in\CONE\colon\,
\sup_{\mu \geq 1}\,\{\phi^{\cI}(\mu)- G(\mu,\overline{x})\}> 0\Big\}.
\end{aligned}
\ee
\ec

As asserted in Theorem~\ref{phtrcurve2} below, if we let $(\alpha,\beta)$ run in $\cD_2$
along a linear segment parallel to the first diagonal, then the free energy $f(\alpha,\beta;p)$
remains constant until $(\alpha,\beta)$ enters $\cL_1$. In other words, if we pick
$(\alpha_0,\beta_0)\in \cD_2$ and consider for $u\geq 0$ the point $s_u=(\alpha_0+u,\beta_0+u)$,
then the free energy $f(s_u;p)$ remains equal to $f(\alpha_0,\beta_0;p)$ until $s_u$ exits $\cD_2$
and enters $\cL_1$. This passage from $\cD_2$ to $\cL_1$ comes with a non-analyticity of the
free energy. This phase transition is represented by a second critical curve in the phase diagram
(see Fig.\ \ref{fig-subcrit}).

\subsubsection{The $\cL_2$-phase: $AB$-localization, $BA$-localization}
\label{S1.4.4}

The remaining phase is:

\bd{L2def}
For $p<p_c$,
\be{L2char}
\cL_2 = \big\{(\alpha,\beta)\in\CONE\colon\, f_{\cL_1}(\alpha,\beta;p)<f(\alpha,\beta;p)\big\}.
\ee
\ed

Starting from $(\alpha,\beta)\in\cL_1$, we increase $\beta$ until it becomes energetically
advantageous for the copolymer to localize at the interface in the $AB$-blocks as well. This
new phase has both $AB$- and $BA$-localization (see Fig.~\ref{fig-subcrits}). Unfortunately,
we are not able to show non-analyticity at the crossover from $\cL_1$ to $\cL_2$ because,
unlike in $\cD_2$, in $\cL_1$ the free energy is not constant in one particular direction
(and the argument we gave for the non-analyticity at the crossover from $\cD_2$ to $\cL_1$
is not valid here). Consequently, the phase transition between $\cL_1$ and $\cL_2$ is still a
\emph{conjecture} at this stage, but we strongly believe that a third critical curve indeed
exists.

\subsection{Main results for the phase diagram}
\label{S1.5}

In Section~\ref{S1.4} we defined the four phases and obtained a characterization of them
in terms of the block pair free energies and the single interface free energy at certain
values of the maximizers in the associated variational formulas. The latter serve as the
starting point for the analysis of the properties of the critical curves (Section~\ref{S1.5.1})
and the phases (Sections~\ref{S1.5.2}--\ref{S1.5.3}).

%%%%%%%%%%%%%%%%%%%%%%%%%%%%%%%%%%%%%%%%%%%%%%%%%%%%%%%%%%%%%%%%%%%%%%%%%%%%%%%%%%%%%%
\begin{figure}[hbtp]
\begin{center}
\setlength{\unitlength}{0.6cm}
\begin{picture}(12,12)(0,-3.5)
\put(0,0){\line(12,0){12}}
\put(0,0){\line(0,8){8}}
\put(0,0){\line(0,-3){3}}
\put(0,0){\line(-3,0){3}}
{\thicklines
\qbezier(0,0)(2,2)(4,4)
\qbezier(4,4)(3.5,2)(3,1)
\qbezier(3,1)(2,0)(1,-1)
\qbezier(3,1)(5,2)(9,3)
\qbezier(4,4)(6,5)(10,6)
}
\qbezier[5](3,0)(3,0.5)(3,1)
\qbezier[20](4,0)(4,2)(4,4)
\qbezier[60](4,4)(6,6)(8,8)
\qbezier[40](0,0)(1.6,-1.5)(3.2,-3.0)
\qbezier[80](0,6.5)(4.5,6.5)(9,6.5)
\put(-.8,-.8){$0$}
\put(12.5,-0.2){$\alpha$}
\put(-0.1,8.5){$\beta$}
\put(0.3,-2.9){$\alpha^*(p)$}
\put(2.5,-2.1){$\alpha^{**}(p)$}
\put(5.9,-1.6){$\alpha_*$}
\qbezier[10](4,0)(4.7,-0.5)(5.4,-1)
\qbezier[8](3,0)(3,-0.65)(3,-1.3)
\qbezier[14](2,0)(1.6,-1)(1.2,-2)
\put(2,0){\circle*{.20}}
\put(4,0){\circle*{.20}}
\put(3,0){\circle*{.20}}
\put(4,4){\circle*{.25}}
\put(3,1){\circle*{.25}}
\put(1.6,.8){$\cD_1$}
\put(6.5,1){$\cD_2$}
\put(6.5,3.5){$\cL_1$}
\put(8.5,7){$\cL_2$}
\put(1.07,-2.3){\vector(-1,-2){0}}
\put(5.7,-1.27){\vector(1,-1){0}}
\put(3,-1.4){\vector(0,-1){0}}
\end{picture}
\end{center}
\caption{\small
Further details of the phase diagram for $p < p_c$ sketched in Fig.\ \ref{fig-subcrit}. There
are four phases, separated by three critical curves, meeting at two tricritical points.}
\label{fig-subcritext}
\end{figure}

%%%%%%%%%%%%%%%%%%%%%%%%%%%%%%%%%%%%%%%%%%%%%%%%%%%%%%%%%%%%%%%%%%%%%%%%%%%%%%%%%%%%

\subsubsection{Critical curves}
\label{S1.5.1}

The first two theorems are dedicated to the critical curves between $\cD_1$ and $\cD_2$,
respectively, between $\cD_2$ and $\cL_1$ (see Fig.~\ref{fig-subcritext}).

\bt{phtrcurve}
Let $p<p_c$.\\
(i) There exists an $\alpha^*(p)\in (0,\infty)$ such that $(\alpha^*(p),0)\in\cD_1$ and
$\cD_1\subset\{(\beta+r,\beta)\colon\; r\leq\alpha^*(p),\beta\geq -\tfrac{r}{2}\}$.\\
(ii) For all $r\in [0,\alpha^{*}(p)]$ there exists a $\beta_c^{1}(r)\geq 0$ such that
$\cD_1\cap\{(\beta+r,\beta)\colon\, \beta\in \R\}$ is the linear segment
\be{linseg1}
\cJ_r^1=\big\{(\beta+r,\beta)\colon\, \beta\in
[-\tfrac{r}{2},\beta_c^{1}(r)]\big\}.
\ee
The free energy $f(\alpha,\beta;p)$ is constant on this segment.\\
(iii) $r \mapsto \beta_c^1(r)$ is continuous on $[0,\alpha^*(p)]$.\\
(iv) Along the curve $r\in(0,\alpha^*(p)]\mapsto (\beta_c^1(r)+r,\beta_c^1(r))$ the two
phases $\cD_1$ and $\cL_1$ touch each other, i.e., for all $r\in (0,\alpha^{*}(p)]$
there exists a $v_r>0$ such that
\be{linseg}
\{(\beta+r,\beta)\colon\, \beta\in (\beta_c^{1}(r),\beta_c^1(r)+v_r]\}\subset\cL_1.
\ee
(v) $\beta^1_c(r)\geq \log(1+(1-e^{-r})^{1/2})$ for all $r\in [0,\alpha^{*}(p)]$.
\et

\bt{phtrcurve2}
Let $p<p_c$.\\
(i) For all $r\in (\alpha^{*}(p),\infty)$ there exists a $\beta_c^{2}(r)>0$ such that
$\cD_2\cap\{(\beta+r,\beta)\colon\, \beta\in \R\}$ is the linear segment
\be{linseg2}
\cJ^2_r=\big\{(\beta+r,\beta)\colon\, \beta\in
[-\tfrac{r}{2},\beta_c^{2}(r)]\big\}.
\ee
The free energy $f(\alpha,\beta;p)$ is constant on this segment.\\
(ii) $r\mapsto \beta_c^2(r)$ is lower semi-continuous on $(\alpha^*(p),\infty)$.\\
(iii) At $\alpha^*(p)$ the following inequality holds:
\be{bc2ineq}
\limsup_{r\downarrow \alpha^{*}(p)} \beta_c^{2}(r)\leq \beta_c^{1}(\alpha^*(p)).
\ee
(iv) There exists an $r_2>\alpha^*(p)$ such that along the interval $(\alpha^*(p),r_2]$
the two phases $\cD_2$ and $\cL_1$ touch each other, i.e., for all $r\in(\alpha^*(p),r_2]$
there exists a $v_r>0$ such that
\be{linseg3}
\{(\beta+r,\beta)\colon\, \beta\in
[\beta_c^{2}(r),\beta_c^2(r)+v_r]\}\subset\cL_1.
\ee
(v) $\beta_c^2(r)\geq \log(1+(1-e^{-r})^{1/2})$ for all $r\in(\alpha^{*}(p),\infty)$.
\et

In \cite{dHW06} it was suggested that the tricritical point where $\cD_1$, $\cD_2$ and
$\cL_1$ meet lies on the horizontal axis. Thanks to Theorem \ref{phtrcurve2}(iii) and (v)
we now know that it lies strictly above.

\subsubsection{Infinite differentiability of the free energy}
\label{S1.5.2}

It was shown in \cite{dHW06}, Lemma 2.5.1 and Proposition 4.2.2, that $f$ is analytic
on the interior of $\cD_1$. We complement this result with the following.

\bt{D2infdiff}
Let $p<p_c$. Then, under Assumption {\rm \ref{assum}} in Section {\rm \ref{S5.3}},
$(\alpha,\beta)\mapsto f(\alpha,\beta;p)$ is infinitely differentiable on the interior of $\cD_2$.
\et

\noindent
Consequently, there are no phase transitions of finite order in the interior of $\cD_1$ and
$\cD_2$.

Assumption \ref{assum} in Section \ref{S5.3} concerns the first supremum in (\ref{fevar}) when
$(\alpha,\beta)\in\cD_2$. Namely, it requires that this supremum is uniquely taken at $(\rho_{kl})
= (\rho^*_{kl}(p))$ with $\rho^*_{AA}(p)+\rho^*_{AB}(p)=\rho^*(p)$ given by (\ref{rho*def}) and
with $\rho^*_{BA}(p)$ maximal subject to the latter equality. In view of Fig.\ \ref{fig-subcrits},
this is a resonable assumption indeed, because in $\cD_2$ the copolymer will first try to maximize
the fraction of time it spends crossing $A$-blocks, and then try to maximize the fraction of time
it spends crossing $B$-blocks that have an $A$-block as neighbor.

We do not have a similar result for the interior of $\cL_1$ and $\cL_2$, simply because we
have insufficient control of the free energy in these regions. Indeed, whereas the variational
formulas (\ref{fevarred}) and (\ref{fevarredD2}) only involve the block free energies $\psi_{AA}$,
$\psi_{BB}$ and $\psi_{BA}^{\widehat\kappa}$, for which (\ref{psiAAandBB}) and (\ref{psiinflink2})
provide closed form expressions, the variational formula in (\ref{fevarredL1}) also involves
the block free energy $\psi_{BA}$, for which no closed form expression is known because
(\ref{psiinflink}) contains the single flat infinite interface free energy $\phi^\cI$.

\subsubsection{Order of the phase transitions}
\label{S1.5.3}

Theorem~\ref{phtrcurve}(ii) states that, in $\cD_1$, for all $r \in [0,\alpha^*(p)]$
the free energy $f$ is constant on the linear segment $\cJ^1_r$, while
Theorem~\ref{phtrcurve2}(i) states that, in $\cD_2$, for all $r \in (\alpha^*(p),\infty)$
the free energy $f$ is constant on the linear segment $\cJ_r^2$. We denote these
constants by $f_{\cD_1}(r)$, respectively, $f_{\cD_2}(r)$.

According to Theorems~\ref{phtrcurve}(ii) and \ref{phtrcurve2}(ii), the phase transition
between $\cD_1$ and $\cD_2$ occurs along the linear segment $\cJ^1_{\alpha^*(p)}$ with
$\beta_c^1(\alpha^*(p))=\alpha^{**}(p)-\alpha^*(p)$. This transition is of order smaller
than or equal to $2$.

\bt{th:orrder}
There exists a $c>0$ such that, for $\delta>0$ small enough,
\be{eq:orderstrait}
c\, \delta^2\leq f_{\cD_2}(\alpha^*(p)+\delta)-f_{\cD_1}(\alpha^*(p))
-f'_{\cD_1}(\alpha^*(p))\,\delta -\frac{1}{2}\,f''_{\cD_1}(\alpha^*(p))\,\delta^2.
\ee
\et

According to Theorem~\ref{phtrcurve}(iv), the phase transition between $\cD_1$ and $\cL_1$
occurs along the curve $\{(r+\beta_c^1(r),\beta_c^1(r))\colon\,r\in[0,\alpha^*(p)]\}$.
This transition is of order smaller than or equal to $2$ and strictly larger than 1.

\bt{th:orrder2}
For all $r \in [0,\alpha^*(p))$ there exist $c>0$ and $\zeta\colon\,[0,1]\mapsto [0,\infty)$
satisfying $\lim_{x\downarrow 0} \zeta(x)=0$ such that, for $\delta>0$ small enough,
\be{eq:order2}
c\, \delta^2\leq f_{\cL_1}(r+\beta_c^1(r)+\delta,\beta_c^1(r)+\delta)
-f_{\cD_1}(r)\leq \zeta(\delta) \delta.
\ee
\et

According to Theorem~\ref{phtrcurve2}(iv), the phase transition between $\cD_2$ and $\cL_1$
occurs at least along the curve
\be{curvetr}
\big\{(r+\beta_c^2(r),\beta_c^2(r))\colon\,r \in [\alpha^*(p),\alpha^*(p)+r_2]\big\}.
\ee
We are not able to determine the exact order of this phase transition, but we can prove
that it is smaller than or equal to the order of the phase transition in the single interface
model. The latter model was investigated (for a different but analogous Hamiltonian) in Giacomin
and Toninelli~\cite{GT06a}, where it was proved that the phase transition is at least of
second order. Numerical simulations suggest that the order is in fact higher than second
order. In what follows we denote by $\gamma$ the order of the single interface transition.
This means that there exist $c_2>c_1>0$ and a slowly varying function $L$ such that, for
$\delta>0$ small enough,
\be{ord}
c_1 \delta^\gamma L(\delta) \leq \phi^{\cI}(\tfrac{c_r}{b_r};r+\beta_c^2(r)
+\delta,\beta_c^2(r)+\delta)-\hat{\kappa}(\tfrac{c_r}{b_r})\leq c_2 \delta^\gamma L(\delta),
\ee
where $(c_r,b_r)$ are the unique maximizers of \eqref{psiinflink2} at $(r+\beta_c^2(r),
\beta_c^2(r);y_r)$ and $y_r$ is the second component of the unique maximizers of \eqref{fevarredD2}
at $(r+\beta_c^2(r),\beta_c^2(r))$.
\bt{th:orrder3}
For all $r\in[\alpha^*(p),\alpha^*(p)+r_2)$ there exist $c>0$ such that, for $\delta>0$ small enough,
\be{eq:order3}
c\, \delta^\gamma L(\delta) \leq f_{\cL_1}(r+\beta_c^2(r)+\delta,\beta_c^2(r)+\delta)
-f_{\cD_2}(r).%\leq c_2\, \delta^\gamma L(\delta).
\ee
\et
We believe that the order of the phase transition along the critical curve separating $\cD_1$ and
$\cD_2$, $\cD_1$ and $\cL_1$, and $\cD_2$ and $\cL_1$ are, respectively, 2, 2 and $\gamma$. However,
except for Theorem \ref{th:orrder2}, in which we give a partial upper bound, we have not been able
to prove upper bounds in Theorems \ref{th:orrder} and \ref{th:orrder3} due to a technical difficulty
associated with the uniqueness of the maximizer $(a_{kl})$ in (\ref{fevar}).

\subsection{Open problems}
\label{S1.6}

The following problems are interesting to pursue (see Fig.~\ref{fig-subcritext}):

\begin{itemize}
\item[(a)]
Prove that  $r\mapsto \beta_c^2(r)$ is continuous on $(\alpha^*(p),\infty)$. Prove that
$r\mapsto \beta_c^1(r)$ is strictly decreasing and $r\mapsto \beta_c^2(r)$ is strictly
increasing.
\item[(b)]
Show that the critical curve between $\cD_2$ and $\cL_1$ meets the critical curve
between $\cD_1$ and $\cD_2$ at the end of the linear segment, i.e., show that
(\ref{bc2ineq}) can be strengthened to an equality.
\item[(c)]
Establish the existence of the critical curve between $\cL_1$ and $\cL_2$.
Prove that the free energy is infinitely differentiable on the interior of
$\cL_1$ and $\cL_2$.
\item[(d)]
Show that the critical curve between $\cD_2$ and $\cL_1$ never crosses the critical
curve between $\cL_1$ and $\cL_2$.
\item[(e)]
Show that the phase transitions between $\cD_1$ and $\cL_1$ and between $\cD_1$ and $\cD_2$
are of order $2$.
\end{itemize}

\subsection{Outline}
\label{S1.7}

In Section~\ref{S2} we derive some preparatory results concerning existence and
uniqueness of maximizers and inequalities between free energies. These will be used
in Section~\ref{S3}--\ref{S4} to prove the claims made in Section~\ref{S1.4}--\ref{S1.5},
respectively.

The present paper concludes the analysis of the phase diagram started in \cite{dHW06}
and continued in \cite{dHP07b}. The results were announced in \cite{dHP07a} without proof.

%%%%%%%%%%%%%%%%%%% SECTION 2 %%%%%%%%%%%%%%%%%%%%%%%%%%%%%%%%%%%%%%%%%%%%%%%%%%%%%%%%%

\section{Preparations}
\label{S2}

\subsection{smoothness of $\hat{\kappa}$ and  $\kappa$}

In this section, we recall some results from \cite{dHW06} concerning the entropies
$\kappa$ and $\hat{\kappa}$ defined in \eqref{kappa} and \eqref{entpath12}.

\bl{smoothkappa} {\rm (\cite{dHW06}, Lemma 2.1.2 and 2.1.1)}\\
(i) $(a,b) \mapsto a\kappa(a,b)$ is continuous and strictly concave
on $\DOM$ and analytic on the interior of $\DOM$.\\
(ii) $\mu\mapsto\mu\hat\kappa(\mu)$ is continuous and strictly concave
on $[1,\infty)$ and analytic on $(1,\infty)$.
\el
This allows to state the following properties of $\psi_{kl}$.
\bc{concpsi1}
(i) For $kl\in \{AA,BB\}$,  $(\alpha,\beta,a)\mapsto\psi_{kl}(\alpha,\beta;a)$ is infinitely
differentiable on $\mathbb{R}^2\times (2,\infty)$.\\
(ii) For $kl\in \{AA,BB\}$ and $(\alpha,\beta)\in\CONE$,  $a\mapsto\psi_{kl}(\alpha,\beta;a)$
is strictly concave on $[2,\infty)$.\\
(iii) For $(\alpha,\beta)\in\CONE$, $a\mapsto \psi_{BA}^{\hat{\kappa}}(\alpha,\beta;a)$ is
strictly concave $[2,\infty)$.
\ec
\bpr
Lemma \ref{smoothkappa} and formulas \ref{psiAAandBB} imply immediately (i) and (ii).
Lemma \ref{smoothkappa} implies also that for all $a\geq 2$, $(c,b)\mapsto c \kappa (c/b)$
and $(c,b)\mapsto (a-c)\kappa(a-c,1-b)$ are strictly concave. The latter, together with
formula \eqref{psiinflink2} are sufficient to obtain (iii).
\epr

\subsection{Smoothness of $\phi^{\cI}$ and $\psi_{kl}$}
\label{S2.4}

In this section, we recall from \cite{dHP07b} some key properties concerning the single
interface free energy and the block pair free energies.

\bl{smo}
(i) $(\alpha,\beta,\mu)\mapsto \phi^{\cI}(\alpha,\beta;\mu)$ is continuous on $\CONE
\times [1,\infty)$.\\
(ii) For all $k,l\in\{A,B\}$, $(\alpha,\beta;a)\mapsto \psi_{kl}(\alpha,\beta;a)$
is continuous on $\CONE \times [2,\infty)$.
\el

\bpr
To prove (i) it suffices to check that $\mu\mapsto \phi^{\cI}(\alpha,\beta;\mu)$ is continuous
on $[1,\infty)$ and that there exists a $K>0$ such that $(\alpha,\beta)\mapsto \phi^{\cI}(\alpha,
\beta;\mu)$ is $K$-Lipshitz for all $\mu\in[1,\infty)$. These two properties are obtained by using,
respectively, the concavity of $\mu\mapsto \mu \phi^{\cI}(\alpha,\beta;\mu)$ and the expression
of the Hamiltonian in \eqref{Zinf}. The proof of (ii) is the same.
\epr

\noindent
Other important results, proven in \cite{dHP07b}, are stated below. They concern the asymptotic
behavior of $\psi_{kl}$, $\phi^{\cI}$ and some of their partial derivatives as $\mu$ and $a$ tend
to $\infty$.

\bl{mulim} {\rm(\cite{dHP07b}, Lemma 2.4.1)}
For any $\beta_0>0$, uniformly in $\alpha\geq \beta$ and $\beta\leq\beta_0$,\\
(i) $\lim_{\mu\to\infty} \phi^\cI(\alpha,\beta;\mu)=0$,\\
(ii) for $kl\in\{AB,BA\}$, $\lim_{a\to\infty} \psi_{kl}(\alpha,\beta;a)= 0$.
\el

\bl{l:c'} {\rm(\cite{dHP07b}, Lemma 5.4.3)}
Fix $(\alpha,\beta) \in \CONE$.\\
(i) For all $k,l\in\{A,B\}$ with $kl\neq BB$, $\lim_{a\to\infty} a\psi_{kl}(a)=\infty$.\\
(ii) Let $\mathcal{K}$ be a bounded subset of $\CONE$. For all $k,l\in\{A,B\}$, $\lim_{a\to\infty}
\partial [a\psi_{kl}(\alpha,\beta;a)]/\partial a\leq 0$ uniformly in $(\alpha,\beta)\in\mathcal{K}$.
\el

\bpr
Only the uniformity in $(\alpha,\beta)\in\mathcal{K}$ in (ii) was not proven in \cite{dHP07b}. This
is obtained as follows. Let $m$ be the minimum of $2\psi_{kl}(\alpha,\beta;2)$ on $\mathcal{K}$. By
Lemma \ref{mulim}(ii), for all $\gep>0$ there exists an $a_0\geq 2$ such that $\psi_{kl}(\alpha,\beta;a)
\leq\gep$ for all $(\alpha,\beta)\in\mathcal{K}$ and $a\geq a_0$. Moreover, by concavity, the derivative
of $a\mapsto a\psi_{kl}(\alpha,\beta;a)$ is decreasing and, consequently, $a\gep-m \geq (a-2)\partial_a
\psi_{kl}(\alpha,\beta;a)$ for $a\geq a_0$. This implies that
\be{uu1}
\partial_a \psi_{kl}(\alpha,\beta;a) \leq
\frac{a\gep-m}{a-2}=\frac{\gep-m/a}{1-2/a}, \qquad a\geq a_0.
\ee
\epr

\subsection{Maximizers for the free energy: existence and uniqueness}
\label{S2.1}

Up to now we have stated the existence and uniqueness of the maximizers of the
variational formula (\ref{fevar}) only in some particular cases. In $\cD_1$ we
recalled the result of \cite{dHW06}, stating the uniqueness of the maximizers $(\overline{x},
\overline{y})$ in the variational formula (\ref{fevarred}), while in $\cD_2$ we
announced the uniqueness of the maximizers $(\overline{x},\overline{y},\overline{z})$
in the variational formula (\ref{fevarredD2}).

For $(\alpha,\beta)\in\CONE$, $p\in(0,1)$ and $ (\rho_{kl}) \in \cR(p)$, let
(recall (\ref{fctV}))
\be{maxdefs}
\begin{aligned}
f_{\alpha,\beta,(\rho_{kl})}
&= \sup_{(a_{kl})\in \cA} V\big((\rho_{kl}),(a_{kl})\big),\\
\cO_{(\rho_{kl})} &= \big\{kl\in\{A,B\}^2\colon\,\rho_{kl}>0\big\},\\
\cJ_{\alpha,\beta,(\rho_{kl})} &= \{(a_{kl})_{kl\in \cO_\rho} \in \cA\colon\,
f_{\alpha,\beta,(\rho_{kl})} = V\big((\rho_{kl}),(a_{kl})\big)\},\\
\cR_{\alpha,\beta,p}^f &= \{(\rho_{kl})\in\cR(p)\colon\,f(\alpha,\beta;p)
= f_{\alpha,\beta,(\rho_{kl})}\},\\
\cP_{\alpha,\beta,p} &= \bigcup_{(\rho_{kl}) \in \cR_{\alpha,\beta,p}^f} \cO_{(\rho_{kl})}.
\end{aligned}
\ee

\bl{p:unimaw}
For every $(\alpha,\beta)\in\CONE$, $p\in (0,1)$ and $(\rho_{kl})\in\cR(p)$, the set
$\cJ_{\alpha,\beta,(\rho_{kl})}$ is non-empty. Moreover, for all $kl\in\cO_{(\rho_{kl})}$
such that $x\mapsto x\psi_{kl}(x)$ is strictly concave, there exists a unique
$a^{(\rho_{kl})}_{kl}\geq 2$ such that $a_{kl}=a^{(\rho_{kl})}_{kl}$ for all $(a_{kl})
\in \cJ_{\alpha,\beta,(\rho_{kl})}$.
\el

\bpr
The proof that $\cJ_{\alpha,\beta,(\rho_{kl})}\neq\emptyset$ is given in \cite{dHP07b},
Proposition 5.5.1. If $(a_{kl})\in\cJ_{\alpha,\beta,(\rho_{kl})}$, then differentiation gives
\be{max1}
\frac{\partial [x\psi_{kl}(x)]}{\partial x}(a_{kl})=f_{\alpha,\beta,(\rho_{kl})},
\ee
which implies the uniqueness of $a_{kl}$ as soon as $x\mapsto x\psi_{kl}(x)$
is strictly concave.
\epr

\br{strconcav2}
{\rm
Note that \eqref{max1} ought really to be written as
\be{partwr}
\partial^- [x\psi_{kl}(x)](a_{kl})\leq f_{\alpha,\beta,(\rho_{kl})}
\leq \partial^+ [x\psi_{kl}(x)](a_{kl}),
\ee
where $\partial^-$ and $\partial^+$ denote the left- and right-derivative. Indeed,
for $kl\in\{AB,BA\}$ we do not know whether $x\mapsto x\psi_{kl}(x)$ is differentiable
or not. However, we know that these functions are concave, which is sufficient to ensure
the existence of the left- and right-derivative. We will continue this abuse of notation
in what follows.
}
\end{remark}

\bp{p:unimaw2}
For every $(\alpha,\beta)\in\CONE$ and $p\in (0,1)$, the set $\cR^f_{\alpha,\beta,p}$ is
non-empty. Moreover, for all  $kl\in \cP_{\alpha,\beta,p}$ such that $x \mapsto x
\psi_{kl}(x;\alpha,\beta)$ is strictly concave, there exists a unique $a_{kl}(\alpha,\beta)
\geq 2$ such that $a^{(\rho)}_{kl}=a_{kl}(\alpha,\beta)$ for all $(\rho_{kl}) \in
\cR_{\alpha,\beta,p}^f$.
\ep

\bpr
We begin with the proof of $\cR^f_{\alpha,\beta,p}\neq \emptyset$. Let $(\rho^B)=(\rho^B_{kl})$
denote the $2\times 2$ matrix with $\rho^B_{BB}=1$ and $\rho^B_{BA}=\rho^B_{AB}=\rho^B_{AA}=0$.

\medskip\noindent
\underline{Case 1}: $\sup_{x\geq 2} \psi_{BB}(x)>0$.\\
Since $\cR(p)$ is a compact set, the continuity of $(\rho_{kl})\mapsto f_{\alpha,\beta,
(\rho_{kl})}$ implies that $\cR^f_{\alpha,\beta,p}\neq\emptyset$. To prove this continuity,
we note that, since $\psi_{kl}\geq \psi_{BB}$ for all $k,l\in\{A,B\}$, $f_{\alpha,\beta,(\rho_{kl})}$
is bounded from below by $\sup_{x\geq 2} \psi_{BB}(x)>0$ uniformly in $(\rho_{kl})\in\cR(p)$.
This is sufficient to mimick the proof of \cite{dHP07b}, Proposition 5.5.1(i), which shows
that there exists a $R>0$ such that, for all $(\rho_{kl})\in \cR(p)$,
\be{falbelb}
f_{\alpha,\beta,(\rho_{kl})}=\sup_{\{(a_{kl})\colon\; a_{kl}\in[2,R]\}}
V((\rho_{kl}),(a_{kl})).
\ee
This in turn is sufficient to obtain the continuity of $(\rho_{kl})\mapsto f_{\alpha,\beta,(\rho_{kl})}$.

\medskip\noindent
\underline{Case 2}: $\sup_{x\geq 2} \psi_{BB}(x)\leq 0$.\\
Since $p>0$ by assumption, we can exclude the case $\cR(p)=\{\rho^B\}$, and therefore we may
assume that $\cR(p)$ contains at least one element different from $(\rho^B)$. Clearly,
$f_{\alpha,\beta,(\rho^B)}\leq 0$, and for any sequence $((\rho_{n}))_{n \geq 1}$ in $\cR(p)$
that converges to $(\rho^B)$ it can be shown that
\be{limsupp}
\limsup_{n\to \infty} f_{\alpha,\beta,(\rho_n)}\leq 0.
\ee
As asserted in Lemma \ref{l:c'}(i), for $kl\neq BB$ we have $\lim_{x\to \infty} x\,\psi_{kl}(x)
=\infty$ and this, together with (\ref{fctV}--\ref{fevar}), forces $f(\alpha,\beta)>0$. Therefore,
\eqref{limsupp} is sufficient to assert that there exists an open neighborhood $\cW$ of $(\rho^B)$
such that $f_{\alpha,\beta,(\rho_{kl})}\leq f(\alpha,\beta;p)/2$ when $(\rho_{kl})\in \cW$, and then
\be{falbesup}
f(\alpha,\beta;p)=\sup_{(\rho_{kl})\in\cR(p)\cup \cW^c} f_{\alpha,\beta,(\rho_{kl})}.
\ee
Finally, $f_{\alpha,\beta,(\rho_{kl})}$ is bounded from below by a strictly positive constant
uniformly in $(\rho_{kl})\in \cR(p)\cup\cW^c$. Hence, by mimicking the proof of Case 1, we obtain that
$(\rho_{kl})\mapsto f_{\alpha,\beta,(\rho_{kl})}$ is continuous on the compact set $\cW^c \cup\cR(p)$.
To complete the proof, we note that, since
\be{impl}
(\rho_1),(\rho_2)\in\cR_{\alpha,\beta,p}^f \quad \Longrightarrow \quad
f_{\alpha,\beta,(\rho_1)}=f_{\alpha,\beta,(\rho_2)},
\ee
(\ref{max1}) implies that $a_{kl}^{(\rho_1)}=a_{kl}^{(\rho_2)}$.
\epr

Proposition \ref{p:unimaw2} gives us the uniqueness of $a_{AA}(\alpha,\beta)$ and $a_{BB}(\alpha,\beta)$
for all $(\alpha,\beta)\in \CONE$. In the following proposition we prove that these functions are
continuous in $(\alpha,\beta)$.

\bp{p:unimaw1}
$(\alpha,\beta)\mapsto a_{AA}(\alpha,\beta)$ and $(\alpha,\beta)\mapsto a_{BB}(\alpha,\beta)$
are continuous on $\CONE$.
\ep

\bpr
Let $kl\in\{AA,BB\}$. By Proposition~\ref{p:unimaw2}, $a_{kl}(\alpha,\beta)$ is the unique solution
of the equation $\partial [x\psi_{kl}(\alpha,\beta;x)]/\partial x=f(\alpha,\beta;p)$. As proved in
Case 2 of Proposition \ref{p:unimaw2}, we have $f(\alpha,\beta;p)>0$. Moreover, with the help
\cite{dHW06}, Lemma 2.2.1, which gives the explicit value of $\kappa(x,1)$, we can easily show that
$\limsup_{x\to \infty} \partial[x\psi_{kl}(\alpha,\beta;x)]/\partial x \leq 0$ uniformly in
$(\alpha,\beta)\in \CONE$. This, together with \eqref{max1} and the fact that $f(\alpha,\beta)$
is bounded when $(\alpha,\beta)$ is bounded, is sufficient to assert that $a_{kl}(\alpha,\beta)$
is bounded in the neighborhood of any $(\alpha,\beta) \in \CONE$. Therefore, by the continuity of
$(\alpha,\beta)\mapsto f(\alpha,\beta)$ and $(\alpha,\beta,x)\mapsto x\psi_{kl}(\alpha,\beta;x)$
and by the uniqueness of $a_{kl}(\alpha,\beta)$ for all $(\alpha,\beta)\in \CONE$, we obtain that
$(\alpha,\beta)\mapsto a_{kl}(\alpha,\beta)$ is continuous.
\epr

\subsection{Inequalities between free energies}
\label{S2.2}

Abbreviate $\cF=\{AA,AB,BA,BB\}$ and let
\be{rhojset}
I = \left\{(\rho_{kl})_{kl\in\cF}\,\colon\,\sum_{kl\in\cF}\rho_{kl}=1,
\,\rho_{kl}>0\,\,\forall\, kl\in\cF\right\}.
\ee
For $kl\in\cF$, let $x\mapsto x\zeta_{kl}(x)$ and $x\mapsto x\overline{\zeta}_{kl}(x)$ be concave
on $[2,\infty)$, $\zeta_{kl}$ be differentiable on $(2,\infty)$, and $\overline{\zeta}_{kl}\geq
\zeta_{kl}$.
For $(\rho_{kl})\in\cR(p)$ and $(x_{kl})\in \cA$, put
\be{vari1}
f_{(\rho_{kl})}((x_{kl}))=\frac{\sum_{kl\in\cF}\,\rho_{kl}\, x_{kl}\,\zeta_{kl}(x_{kl})}
{\sum_{kl\in\cF}\,\rho_{kl}\, x_{kl}}
\quad \text{ and } \quad
\overline{f}_{(\rho_{kl})}((x_{kl}))
= \frac{\sum_{kl\in\cF}\,\rho_{kl}\, x_j\,\overline{\zeta}_{kl}(x_{kl})}{\sum_{kl\in\cF}\,\rho_{kl}\, x_{kl}}
\ee
and
\be{vari}
f=\sup_{(\rho_{kl})\in\cR(p)}
\sup_{(x_{kl})\in\cA}f_{(\rho_{kl})}((x_{kl}))
\quad \text{ and } \quad
\overline{f}=\sup_{(\rho_{kl})\in\cR(p)}
\sup_{(x_{kl})\in \cA} \overline{f}_{(\rho_{kl})}((x_{kl})).
\ee

\bl{lemdu}
Assume that there exist $(\tilde{\rho}_{kl})\in\cR(p)\cap I $ and $(\tilde{x}_{kl})\in (2,\infty)^{\cF}$
that maximize the first variational formula in {\rm (\ref{vari})}. Then the following
are equivalent:\\
(i) $\overline{f}>f$;\\
(ii) there exists a $kl\in\cF$ such that $\overline{\zeta}_{kl}(\tilde{x}_{kl})>\zeta_{kl}(\tilde{x}_{kl})$.
\el

\bpr
This proposition is a generalization of \cite{dHW06}, Proposition 4.2.2. It is obvious that
(ii) implies (i). Therefore it will be enough to prove that $\overline{f}=f$ when (ii) fails.
Trivially, $\overline{f}\geq f$.

Abbreviate $\theta_{kl}(x)=x\zeta_{kl}(x)$ and $\overline{\theta}_{kl}(x)=x\overline{\zeta}_{kl}(x)$.
If (ii) fails, then $\theta_{kl}(\tilde{x}_{kl})=\overline{\theta}_{kl}(\tilde{x}_{kl})$ for all
$kl \in\cF$. Since, by assumption, $\theta_{kl}$ is differentiable, $\theta_{kl}$ and
$\overline{\theta}_{kl}$ are concave and $\overline{\theta}_{kl}\geq\theta_{kl}$, it
follows that $\overline{\theta}_{kl}$ is differentiable at $\tilde{x}_{kl}$ with
$(\overline{\theta}_{kl})'(\tilde{x}_{kl})=(\theta_{kl})'(\tilde{x}_{kl})$. The fact that
$(\tilde{\rho}_{kl})\in\cR(p)\cap I$ and $(\tilde{x}_{kl})\in (2,\infty)^\cF$ maximize the first
variational formula in (\ref{vari1}) implies, by differentiation of the l.h.s.\ of
(\ref{vari1}) w.r.t.\ $\tilde{x}_{kl}$ at $((\tilde{\rho}_{kl}),(\tilde{x}_{kl}))$, that
$(\theta_{kl})'(\tilde{x}_{kl})=f$ for all $kl\in\cF$. Therefore $(\overline{\theta}_{kl})'
(\tilde{x}_{kl})=f$ for all $kl\in\cF$. Now pick $(\rho_{kl})\in\cR(p)$, $(x_{kl})\in\cA$,
and put $N=\sum_{kl\in\cF}\rho_{kl}\theta_{kl}(\tilde{x}_{kl})$, $V=\sum_{kl\in\cF}\rho_{kl}
\tilde{x}_{kl}$.
Since $\overline{\theta}_{kl}$ is concave, we can write
\be{arti}
\begin{aligned}
\overline{f}_{(\rho_{kl})}((x_{kl}))
&=\frac{N+\sum_{kl\in\cF}\,\rho_{kl}\,(\overline{\theta}_{kl}(x_{kl})-\overline{\theta}_{kl}
(\tilde{x}_{kl}))}
{V+\sum_{kl\in\cF}\,\rho_{kl}\, (x_{kl}-\tilde{x}_{kl})}
\leq \frac{N+f \sum_{kl\in\cF}\,\rho_{kl}\, (\tilde{x}_{kl}-x_{kl})}
{V+\sum_{kl\in\cF}\,\rho_{kl}\, (x_{kl}-\tilde{x}_{kl})}.
\end{aligned}
\ee
But $N/V=f_{(\rho_{kl})}((\tilde{x}_{kl}))\leq f$, and therefore (\ref{arti}) becomes
$\overline{f}_{(\rho_{kl})}((x_{kl}))\leq f$, which, after taking the supremum over $(\rho_{kl})
\in\cR\cap I$ and $(x_{kl})\in\cA$, gives us $\overline{f}\leq f$.
\epr

%%%%%%%%%%%%% SECTION 3 %%%%%%%%%%%%%%%%%%%%%%%%%%%%%%%%%%%%%%%%%%%%%%%%%%%%%%%%%%%%%%%%%%%%

\section{Characterization of the four phases}
\label{S3}

\subsection{Proof of Proposition \ref{d2char}}
\label{S3.1}

\bpr
Recall that $(\overline{x},\overline{y},\overline{z})$ is the unique maximizer of the variational
formula in (\ref{fevarredD2}) at $(\alpha,\beta;p)$. By (\ref{phaseD1}), $f=f_{\cD_1}$ if
$(\alpha,\beta)\in\cD_1$ and $f>f_{\cD_1}$ otherwise, and therefore Proposition \ref{d2char} will
be proven if we can show that
\be{fpsirels}
\begin{aligned}
&\psi_{AB}(\overline{x})=\psi_{AA}(\overline{x})
\quad \mbox{ and } \quad
\psi_{BA}(\overline{y})=\psi_{BA}^{\hat{\kappa}}(\overline{y})
\quad &\Longrightarrow \quad f=f_{\cD_2},\\
&\psi_{AB}(\overline{x})>\psi_{AA}(\overline{x})
\quad \mbox{ or } \quad
\psi_{BA}(\overline{y})>\psi_{BA}^{\hat{\kappa}}(\overline{y})
\quad &\Longrightarrow \quad f>f_{\cD_2}.
\end{aligned}
\ee
But this follows by applying Lemma \ref{lemdu} with $\zeta_{kl}=\overline{\zeta}_{kl}
=\psi_{kl}$ for $kl\in\{AA,BB\}$, $\zeta_{BB}=\overline{\zeta}_{BB}=\psi_{BB}$,
$\zeta_{AB}=\psi_{AA}$, $\overline{\zeta}_{AB}=\psi_{AB}$, $\zeta_{BA}=\psi_{BA}^{\hat{\kappa}}$
and $\overline{\zeta}_{BA}=\psi_{BA}$.
\epr

\subsection{Proof of Corollary \ref{philbextext3}}
\label{S3.2}

\bpr
By Lemma~\ref{lemdu1}, $\psi_{AB}(\overline{x})>\psi_{AA}(\overline{x})$ if and only if
$\sup_{\mu \geq 1}\,\{\phi^{\cI}(\mu)- G(\mu,\overline{x})\}>0$, with $G(\mu,x)$ defined
in (\ref{philb4}). Combine this with Lemma~\ref{linint} below at $\overline{y}$.
\epr

\bl{linint}
For all $y\geq 2$, $\psi_{BA}(y)>\psi^{\hat{\kappa}}_{BA}(y)$ if and only of $\phi^{\cI}
(\tilde{c}/\tilde{b})>\hat{\kappa}(\tilde{c}/\tilde{b})$ with $(\tilde{b},\tilde{c})$
the unique maximizer of the variational formula {\rm (\ref{psiinflink2})} for
$\psi_{BA}^{\hat{\kappa}}(y)$.
\el

\bpr
If $\phi^{\cI}(\tilde{c}/\tilde{b})>\hat{\kappa}(\tilde{c}/\tilde{b})$, then clearly
$\psi_{BA}(y)>\psi^{\hat{\kappa}}_{BA}(y)$. Thus, it suffices to assume that $\psi_{BA}(y)
>\psi^{\hat{\kappa}}_{BA}(y)$ and $\phi^{\cI}(\tilde{c}/\tilde{b})=\hat{\kappa}(\tilde{c}
/\tilde{b})$ and show that this leads to a contradiction. For $(b,c)\in\DOM(y)$,
let
\be{set1}
\begin{aligned}
T(b,c)
&= c \phi^{\cI}(c/b)+(y-c)\big\{\kappa(y-c,1-b)+\tfrac12(\beta-\alpha)\big\},\\
T^{\hat{\kappa}}(b,c)
&= c \hat{\kappa}(c/b)+(y-c)\big\{\kappa(y-c,1-b)+\tfrac12(\beta-\alpha)\big\}.
\end{aligned}
\ee
By definition, the unique maximizer of $T^{\hat{\kappa}}$ on $\DOM(y)$ is $(\tilde{b},
\tilde{c})$. Moreover, $\phi^{\cI}(\tilde{c}/\tilde{b})=\hat{\kappa}(\tilde{c}/\tilde{b})$
implies that $T(\tilde{b},\tilde{c})=T^{\hat{\kappa}}(\tilde{b},\tilde{c})$.
However, $\psi_{BA}(y)>\psi^{\hat{\kappa}}_{BA}(y)$ implies that there exists a $(b',c')
\in\DOM(z)$ such that $T(b',c')>T(\tilde{b},\tilde{c})$. Now
put
\be{zetadefpr}
\zeta\colon\, t \mapsto (\tilde{b},\tilde{c}) + t\,(b'-\tilde{b},c'-\tilde{c}).
\ee
Since $(b,c)\mapsto T^{\hat{\kappa}}(b,c)$ is differentiable and concave on $\DOM(y)$ (recall
that $\hat{\kappa}$ and $\kappa$ are differentiable), also $t\mapsto T^{\hat{\kappa}}(\zeta(t))$
is differentiable and concave, and reaches its maximum at $t=0$. Moreover, $t\mapsto
T(\zeta(t))$ is concave and, since $T(\zeta(\cdot)) \geq T^{\hat{\kappa}}(\zeta(\cdot))$
and  $T(\zeta(0))=T^{\hat{\kappa}}(\zeta(0))$, it follows that $t\mapsto T(\zeta(t))$ is
differentiable at $t=0$ with zero derivative. It therefore is impossible that $T(\zeta(1))
>T(\zeta(0))$.
\epr

\subsection{Proof of Proposition \ref{L1char}}
\label{S3.3}

\bpr
Recall that $(\overline{x},\overline{y},\overline{z})$ is the unique maximizer of the variational
formula in (\ref{fevarredL1}) at $(\alpha,\beta;p)$. By (\ref{phaseD1}) and (\ref{phaseD2}),
$f>f_{\cD_2}$ if $(\alpha,\beta)\in\cD_1^c\cap\cD_2^c$ and $f=f_{\cD_2}$ otherwise. To prove
Proposition \ref{L1char}, we must show that
\be{fpsirelsalt}
\begin{aligned}
&\psi_{AB}(\overline{x})=\psi_{AA}(\overline{x})
\quad &\Longrightarrow \quad f=f_{\cL_1},\\
&\psi_{AB}(\overline{x})>\psi_{AA}(\overline{x})
\quad &\Longrightarrow \quad f>f_{\cL_1}.
\end{aligned}
\ee
But this follows by applying Lemma~\ref{lemdu} with $\zeta_{kl}=\overline{\zeta}_{kl}=\psi_{kl}$
for $kl\in\{AA,BB,BA\}$, $\zeta_{AB}=\psi_{AA}$ and $\overline{\zeta}_{AB}=\psi_{AB}$.
\epr

\subsection{Proof of Corollary \ref{corola3}}
\label{S3.4}

\bpr
This follows by applying Lemma~\ref{lemdu1} to $\psi_{AB}(\overline{x})$.
\epr

%%%%%%%%%%%%%%% SECTION 4 %%%%%%%%%%%%%%%%%%%%%%%%%%%%%%%%%%%%%%%%%%%%%

\section{Proof of the main results for the phase diagram}
\label{S4}

\subsection{Proof of Theorem \ref{phtrcurve}}
\label{S4.1}

In what follows, we abbreviate $\alpha^*=\alpha^*(p)$ and $\rho^*=\rho^*(p)$. We
recall the following.

\bp{p:Fanal} {\rm (\cite{dHW06}, Proposition 2.5.1)}\\
Let $(\alpha,\beta)\in\CONE$ and $\rho\in(0,1)$. Abbreviate $C=\alpha-\beta\geq 0$.
The variational formula in {\rm (\ref{fevarred})} has unique maximizers $\bar x
=\bar x(C,\rho)$ and $\bar y=\bar y(C,\rho)$ satisfying:\\
(i) $2 < \bar y < a^* < \bar x < \infty$ when $C>0$ and $\bar x=\bar y=a^*$
when $C=0$.\\
(ii) $u(\bar x)>v(\bar y)$ when $C>0$ and $u(\bar x) = v(\bar y)$ when $C=0$.\\
(iii) $\rho\mapsto\bar x(C,\rho)$ and $\rho\mapsto\bar y(C,\rho)$ are analytic
and strictly decreasing on $(0,1)$ for all $C>0$.\\
(iv) $C\mapsto\bar x(C,\rho)$ and $C\mapsto\bar y(C,\rho)$ are analytic and
strictly increasing, respectively, strictly decreasing on $(0,\infty)$
for all $\rho\in(0,1)$.
\ep

We are now ready to give the proof of Theorem~\ref{phtrcurve}.

\bpr
(i) Let $(\overline{x},\overline{y})$ be the maximizer of the variational formula in
(\ref{fevarred}) at $(\alpha,\beta)$. Recall the criterion (\ref{philbextext}), i.e.,
\be{kappacrit}
(\alpha,\beta)\in \cD_1^c\quad \text{ if and only if } \quad
\sup_{\mu \geq 1} \Big\{\phi^{\cI}(\alpha,\beta;\mu) +
\tfrac12(\alpha-\beta) - G(\mu,\overline{y})\Big\}>0.
\ee
Since $\phi^{\cI}(\alpha,0;\mu)=\hat{\kappa}(\mu)$ for all $\alpha\geq 0$ and $\mu\geq 1$,
the r.h.s.\ in (\ref{kappacrit}) can be replaced, when $\beta=0$, by
\be{kappacrit2}
\sup_{\mu \geq 1} \Big\{\hat{\kappa}(\mu)+ \frac12\alpha- G(\mu,\overline{y})\Big\}>0.
\ee
Since, by Proposition \ref{p:Fanal}, $\overline{y}$ depends on $C=\alpha-\beta$ only,
the same is true for the l.h.s.\ in (\ref{kappacrit2}). Moreover, as shown in \cite{dHW06},
Proposition 4.2.3(iii), the l.h.s.\ of (\ref{kappacrit2}) is strictly negative at $C=0$,
strictly increasing in $C$ on $[0,\infty)$, and tends to infinity as $C\to\infty$. Therefore
there exists an $\alpha^*\in (0,\infty)$ such that the l.h.s.\ in (\ref{kappacrit2})
is strictly positive if and only if $\alpha-\beta>\alpha^*$. This implies that
$(\alpha^*,0)\in\cD_1$ and, since $\phi^{\cI}(\alpha,\beta;\mu)\geq \hat{\kappa}(\mu)$
for all $(\alpha,\beta)\in \R^2$ and $\mu \geq 1$, it also implies that $(\beta+r,\beta)
\in\cD_1^c$ when $r>\alpha^*$ and $\beta\geq -\frac{r}{2}$.

\medskip\noindent
(ii) The existence of $\beta_c^1(r)$ is proven in \cite{dHW06}, Theorem 1.5.3(ii).
Consequently, the segment $\cJ^1_r=\{(\beta+r,\beta)\colon\,\beta\in[-\tfrac{r}{2},
\beta_c^{1}(r)]\}$ is included in $\cD_1$. This means that $f(\alpha,\beta;p)$ is constant
and equal to $f_{\cD_1}(r)$ on $\cJ^1_r$.

\medskip\noindent
(iii) The continuity of $r\mapsto \beta_c^1(r)$ is proven in \cite{dHW06}, Theorem 1.5.3(ii).

\medskip\noindent
(iv)
Let $r\leq\alpha^*$ and, for $u>0$ let $s_u=(r+\beta_c^1(r)+u,\beta_c^1(r)+u)$. By the
definition of $\beta_c^1(r)$, we know that $s_u\in\cD_1^c$ and therefore that $f(s_u;p)
>f_{\cD_1}(s_u)$. Moreover, since $f_{\cD_2}$ depends only on $\alpha-\beta$, $f(s_u;p)$
cannot be equal to $f_{\cD_2}(s_\mu;p)$, otherwise $f$ would be constant on
$\cJ^1_{\beta^1_c(r)+u}$ (which would contradict the definition of $\beta_c^1(r)$).
Thus, denoting by $(\overline{x}_u,\overline{y}_u,\overline{z}_u)$ any maximizer of
(\ref{fevarredL1}) at $s_u$ (recall that $\overline{x}_u$ is unique by Proposition
\ref{p:unimaw2}), if we prove that there exists a $v>0$ such that $\psi_{AB}(s_u,\overline{x}_u)
=\psi_{AA}(s_u,\overline{x}_u)$ when $u\in(0,v]$, then Proposition \ref{L1char} implies
that $f(s_u;p)=f_{\cL_1}(s_u;p)$. Since $s_0\in\cD_1$, we know from \cite{dHW06}, Proposition
4.2.3(i), that
\be{supe}
\sup_{\mu\geq 1} \big\{\phi^{\cI}(\mu;s_0)-G(\mu,\overline{x}_0)\big\} < 0.
\ee
It follows from \cite{dHP07b}, Lemma 2.4.1, that $\phi^{\cI}(\mu;s_u)\to 0$ as $\mu\to\infty$
uniformly in $u \in [0,1]$ on the linear segment $\{s_u\colon\,u\in[0,1]\}$. Moreover,
Proposition~\ref{p:unimaw1} implies that $u\mapsto \overline{x}_u$ is continuous and
Proposition \ref{p:Fanal}(i) that $\overline{x}_u>a^*=5/2$ for all $u\in [0,1]$. Then,
since $G(\mu,\overline{x}_u)\geq 1/4 \log[\overline{x}_u/(\overline{x}_u-2)]$, we can
assert that there exists an $R>0$ and a $\mu_0>1$ such that $\sup_{\mu \geq \mu_0}\{\phi^{\cI}
(\mu;s_u)-G(\mu,\overline{x}_u)\}\leq -R$ for all $u\in[0,1]$. Moreover, by Lemma \ref{smo}(i)
and by \eqref{philb4} we know that $(\mu,u)\mapsto \phi^{\cI}(\mu;s_u)-G(\mu,\overline{x}_u)$
is continuous and strictly negative on the set $[1,\mu_0]\times\{0\}$. Therefore we can choose
$v>0$ small enough so that $\sup_{\mu\in[1,\mu_0]}\{\phi^{\cI}(\mu;s_u)-G(\mu,\overline{x}_u)\}<0$
for $u\in[0,v]$.

\medskip\noindent
(v) For $r\geq 0$, let $\cT_r=\big\{(\beta+r,\beta)\colon\,\beta\in[-\tfrac{r}{2},\log(1
+\sqrt{1-e^{-r}})]\}$. By an annealed computation we can prove that, for all $r
\geq 0$, $(\alpha,\beta)\in \cT_r$ implies $\phi^{\cI}(\alpha,\beta;\mu)=\hat{\kappa}(\mu)$
for all $\mu\geq 1$. Consequently, the criterion given in Corollary~\ref{philbextext} (for
$(\alpha,\beta) \in \cD_1^c$) reduces to $\sup_{\mu \geq 1}\{\hat{\kappa}(\mu)+ \frac{r}{2}
- G(\mu,\overline{y})\}>0$. By definition of $\alpha^*$, this criterion is not satisfied
when $r\leq \alpha^*$, and therefore $\cT_r\subset \cD_1$. Hence, $\beta_c^{1}(r)\geq
\log(1+\sqrt{1-e^{-r}})$.
\epr

\subsection{Proof of Theorem \ref{phtrcurve2}}
\label{S4.2}

Below we suppress the $p$-dependence of the free energy to ease the notation.

\bpr
(i) From Theorem~\ref{phtrcurve}(i) we know that $f(\beta+r,\beta)>f_{\cD_1}(r)$
when $r>\alpha^*$ and $\beta\geq -\tfrac{r}{2}$. Hence we must show that for all
$r\in(\alpha^*,\infty)$ there exists a $\beta_c^2(r)$ such that $f(\beta+r,\beta)
=f_{\cD_2}(r)$ when $\beta\in[-\frac{r}{2},\beta_c^2(r)]$ and $f(\beta+r,\beta)>f_{\cD_2}(r)$
when $\beta>\beta_c^2(r)$. This is done as follows. Since $\phi^{\cI}(\beta+r,\beta;\mu)
=\hat{\kappa}(\mu)$ for all $\mu\geq 1$ and $-r/2\leq\beta\leq 0$, we have $\psi_{AB}
(\beta+r,\beta;a)=\psi_{AA}(\beta+r,\beta;a)$ and $\psi_{BA}(\beta+r,\beta;a)
=\psi_{BA}^{\hat{\kappa}}(\beta+r,\beta;a)$ for all $a\geq 2$. Therfore Proposition
\ref{d2char} implies $f(\beta+r,\beta)=f_{\cD_2}(r)$ for all $-r/2\leq\beta\leq 0$.
Moreover, $\beta\mapsto f(\beta+r,\beta)$ is convex and therefore the proof will be
complete once we show that there exists a $\beta>0$ such that $f(\beta+r,\beta)>f_{\cD_2}(r)$.
To prove the latter, we recall Corollary~\ref{philbextext3}, which asserts that $(\beta+r,\beta)
\in\cD_2^c$ in particular when
\be{critAB}
\sup_{\mu \geq 1} \Big\{\phi^{\cI}(\mu)-G(\mu,\overline{x})\Big\}>0,
\ee
where $(\overline{x},\overline{y},\overline{z})$ is the maximizer of (\ref{fevarredD2}) at
$(\beta+r,\beta)$, which depends on $r$ only. It was shown in \cite{dHW06}, Equation (4.1.17),
that $\phi^{\cI}(\alpha,\beta;\frac{9}{8})\geq \frac{\beta}{8}$. Therefore, for $r>\alpha^*$
and $\beta$ large enough, the criterion in (\ref{critAB}) is satisfied at $(\beta+r,\beta)$.
Finally, since $f_{\cD_2}$ is a function of $\alpha-\beta$ and $p$ and since $f(\alpha,\beta;p)
=f_{\cD_2}(\alpha-\beta;p)$ for all $(\alpha,\beta)\in\cJ_r^2$, it follows that the free energy
is constant on $\cJ_r^2$.

\medskip\noindent
(ii) To prove that $r\mapsto\beta_c^{2}(r)$ is lower semi-continuous, we must show that for
all $x\in(\alpha^*,\infty)$
\be{limsup}
\limsup_{r \to x} \beta_c^2(r) \leq \beta_c^2(x).
\ee
Set $l=\limsup_ {r\to x} \beta_c^2(r)$. Then there exists a sequence $(r_n)$ with $\lim_{n\to\infty}
r_n=x$ and $\lim_{n\to\infty} \beta_c^2(r_n)=l$. We note that $(\alpha,\beta)\mapsto f(\alpha,\beta)$
and $(\alpha,\beta)\mapsto f_{\cD_2}(\alpha-\beta)$ are both convex and therefore are both continuous.
Effectively, as in (\ref{fedef}), $f_{\cD_2}$ can be written as the free energy associated with the
Hamiltonian in (\ref{Hamiltonian}) and with an appropriate restriction on the set of paths
$\cW_{n,L_n}$, which implies its convexity. By the definition of $\beta_c^2(r_n)$, we can assert
that $f-f_{\cD_2}=0$ on the linear segment $\cJ^2_{r_n}$ for all $n$. Thus, by the continuity of
$(\alpha,\beta)\mapsto (f-f_{\cD_2})(\alpha,\beta)$ and by the convergence of $r_n$ to $x$, we can
assert that $f-f_{\cD_2}=0$ on $\{(\beta+x,\beta)\colon\,-\tfrac{x}{2}\leq\beta\leq l\}$, which
implies that $l\leq \beta_c^2(x)$ by the definition of $\beta_c^2(x)$.

\medskip\noindent
(iii) Set $l=\limsup_{r\to \alpha^*}\beta_c^2(r)$. In the same spirit as the proof of (ii), since
$f-f_{\cD_2}$ is continuous and equal to $0$ on every segment $\cJ^2_r$, it must be that $f$ is
constant and equal to $f_{\cD_2}(\alpha^*)$ on the segment $\{(\beta+\alpha^*,\beta)\colon\,
-\tfrac{\alpha^*}{2}\leq \beta\leq l\}$. This, by the definition of $\beta_c^1(\alpha^*)$, implies
that $l\leq \beta_c^1(\alpha^*)$.

\medskip\noindent
(iv) We will prove that there exist $r_2>\alpha^*$ and $\eta>0$ such that, for all $r\in(\alpha^*,r_2)$
and all $u\in[0,\eta]$,
\be{rala}
\sup_{\mu\geq 1}\big\{\phi^{\cI}(r+\beta_c^2(r)+u,\beta_c^2(r)+u;\mu)
- G(\mu,\overline{x}_{r,u})\big\}\leq 0,
\ee
where $\overline{x}_{r,u}$ is the first coordinate of the maximizer of \eqref{fevarredL1}
at $(r+\beta_c^2(r)+u,\beta_c^2(r)+u)$. This is sufficient to yield the claim, because by
Corollary \ref{corola3} it means that $f_{\cL_1}=f$.

By using (iii), as well as (v) below, we have
\be{linsupinf}
0<\liminf_{r\downarrow \alpha^{*}} \beta_c^{2}(r) \leq
\limsup_{r\downarrow \alpha^{*}} \beta_c^{2}(r) \leq \beta_c^{1}(\alpha^*),
\ee
and hence for all $\gep>0$ there exists a $r_\gep>\alpha^*$ such that, for $\alpha^*<r<r_\gep$,
\be{rala2}
0\leq\beta_c^{2}(r)\leq \beta_c^{1}(\alpha^*)+\gep.
\ee
Next, we define the function
\be{lsuper}
L\colon\,(\alpha,\beta;\mu)\in \CONE\times [1,\infty)
\mapsto \phi^{\cI}(\alpha,\beta;\mu)- G(\mu,\overline{x}_{\alpha,\beta}),
\ee
where $\overline{x}_{\alpha,\beta}$ is the first coordinate of the maximizer of \eqref{fevarredL1}
at $(\alpha,\beta)$, and we set
\be{lsuper2}
F\colon\,(\alpha,\beta)\in \CONE \mapsto \sup_{\mu\geq 1} L(\alpha,\beta;\mu).
\ee
We will show that there exist $r_1>\alpha^*$ and $v>0$ such that $F(\alpha,\beta)$ is non-positive
on the set $\{(r+u,u)\colon\; r\in[\alpha^*,r_1], u\in[0,\beta_c^1(\alpha^*)+v]\}$. Thus, choosing
$\gep=\tfrac{v}{2}$ in \eqref{rala2}, and $r_{2}=\min\{r_{v/2},r_1\}$ and $\eta=\tfrac{v}{2}$ in
\eqref{rala}, we complete the proof.

In what follows we abbreviate $\beta^*=\beta_c^1(\alpha^*)$, $I_1=[\alpha^*,\alpha^*+1]$ and
$I_2=[0,\beta^*+1]$. Since $(\alpha^*+\beta^*,\beta^*)\in\cD_1$, we know from \cite{dHW06},
Proposition 4.2.3(i), that $F(\alpha^*+\beta^*,\beta^*)<0$. Moreover, $\overline{x}_{\alpha^*+u,u}$
is equal to $\overline{x}_{\alpha^*+\beta^*,\beta^*}$ for $u\leq \beta^*$ and, by convexity,
$u\mapsto \phi^{\cI}(\alpha^*+u,u;\mu)$ is non-decreasing for all $\mu\geq 1$. This implies
that $F(\alpha^*+u,u)\leq F(\alpha^*+\beta^*,\beta^*)$ for all $u\in [0,\beta^*]$. Then,
mimicking the proof of Theorem~\ref{phtrcurve}(iv), we use Lemma \ref{mulim}, which tells
us that $\phi^{\cI}(\alpha+\beta,\beta;\mu)\to 0$ as $\mu\to\infty$ uniformly in $(\alpha,\beta)
\in I_1\times I_2$. Moreover, Proposition~\ref{p:unimaw1} implies that $(\alpha,\beta)\mapsto
\overline{x}_{\alpha+\beta,\beta}$ is continuous and, since $G(\mu,x)\geq 1/4 \log[x/(x-2)]$
for $x\geq 2$, we have that there exist $R>0$ and $\mu_0>1$ such that, for all $(\alpha,\beta)
\in I_1\times I_2$,
\be{N11}
\sup_{\mu \geq \mu_0}\{\phi^{\cI}(\alpha+\beta,\beta;\mu)
-G(\mu,\overline{x}_{\alpha+\beta,\beta})\}\leq -R.
\ee
Note that, by \eqref{philb4} and Lemma \ref{smo}(i), the function $L$ defined in \eqref{lsuper}
is continuous on $\CONE\times [1,\infty)$. Moreover, $L(\alpha^*+u,u;\mu)\leq F(\alpha^*+\beta^*,
\beta^*)<0$ for all $\mu \geq 1$ and $u\in [0,\beta^*]$. Therefore, by the continuity of $L$, we
can choose $R_1>0$, $r_2>\alpha^*$ and $v>0$ small enough such that
\be{N12}
\sup_{\mu\leq \mu_0} L(\alpha+\beta,\beta;\mu)\leq -R_1
\ee
for $\alpha\in [\alpha^*,r_2]$ and $\beta\in [0,\beta^*+v]$.

\medskip\noindent
(v) For $r\geq 0$, let
\be{cTRdef}
\cT_r = \Big\{(\beta+r,\beta)\colon\,\beta\in\big[-\tfrac{r}{2},
\log\big(1+(1-e^{-r})^{\frac{1}{2}}\big)\big]\Big\}.
\ee
By an annealed computation we can show that, for all $r\geq 0$, $(\alpha,\beta)\in \cT_r$
implies $\phi^{\cI}(\alpha,\beta;\mu)=\hat{\kappa}(\mu)$ for all $\mu\geq 1$. Moreover,
$\phi^{\cI}\equiv\hat{\kappa}$ implies $\psi_{AB}\equiv\psi_{AA}$. Therefore, for $r>\alpha^*$,
using the criterion (\ref{intersec2}), we obtain that $\cT_r\subset\cD_2$, because none of the
conditions for $(\alpha,\beta)$ to belong to $\cD_1^{c}\cap\cD_2^c$ are satisfied in $\cT_r$.
Hence $\beta_c^{2}(r)\geq \log(1+\sqrt{1-e^{-r}})$.
\epr

\subsection{Proof of Theorem \ref{D2infdiff}}
\label{infdiff}

In this section we give a sketch of the proof of the infinite differentiability of
$(\alpha,\beta)\mapsto f(\alpha,\beta;p)$ on the interior of $\cD_2$. For that, we
mimick the proof of \cite{dHP07b}, Theorem 1.4.3, which states that, in the supercritical
regime $p \geq p_c$, the free energy is infinitely differentiable throughout the localized
phase. The details of the proof are very similar, which is why we omit the details.

It was explained in Section \ref{S1.4.2} that, throughout $\cD_2$, all the quantities involved
in the variational formula in (\ref{fevarredD2}) depend on $(\alpha,\beta)$ only through the
difference $r=\alpha-\beta$. Therefore, it suffices to show that $r\mapsto f_{\cD_2}(r)$
(defined at the beginning of Section \ref{S1.5.3}) is infinitely differentiable on
$(\alpha^{*}(p),\infty)$.

\subsubsection{Smoothness of $\psi_{BA}^{\hat{\kappa}}$ in its localized phase}
\label{S5.3}

This section is the counterpart of \cite{dHP07b}, Section 5.4. Let
\begin{equation}
\cL_{\psi^{\hat{\kappa}}} = \{(r,a)\in (\alpha^*,\infty) \times [2,\infty)\colon\,
\psi_{BA}^{\hat{\kappa}}(r;a)>\psi_{BB}(r;a)\},
\end{equation}
where $\psi_{BB}(r;a)=\kappa(a,1)-\tfrac{r}{2}$ (recall (\ref{psiAAandBB})). Our main result
in this section is the following.

\bp{psiklsmooth}
$(\alpha,\beta,a)\mapsto\psi_{BA}^{\hat{\kappa}}(\alpha,\beta;a)$ is infinitely differentiable
on $\cL_{\psi^{\hat{\kappa}}}$.
\ep

\bpr
Let ${\rm int}[\DOM(a)]$ be the interior of $\DOM(a)$. The proof of the infinite differentiability
of $\psi_{AB}$ on the set
\be{def}
\big\{(\alpha,\beta,a)\in\CONE\times [2,\infty)\colon\;
\psi_{AB}(\alpha,\beta;a)>\tfrac12 \log\tfrac52\big\},
\ee
which was introduced in \cite{dHP07b}, Section 5.4.1, can be readily extended after replacing
$\psi_{AB}$ and $\phi^{\cI}$ on their domains of definition by $\psi_{BA}^{\hat{\kappa}}$
on $\cL_{\psi^{\hat{\kappa}}}$, respectively, $\hat{\kappa}$ on ${\rm int}[\DOM(a)]$. For this
reason, we will only repeat the main steps of the proof and refer to \cite{dHP07b}, Section 5.4.1,
for details.

We begin with some elementary observations. Fix $r \in (\alpha^*,\infty)$, and recall that the supremum
of the variational formula in (\ref{psiinflink}) is attained at a unique pair $(c(r,a),b(r,a))\in
{\rm int}[\DOM(a)]$. Let
\be{partialder}
F(c,b)=c\hat{\kappa}(c/b),\qquad \tilde{F}(c,b)=(a-c)\big[\kappa(a-c,1-b)-\tfrac{r}{2}\big],
\ee
and denote by $F_{c},F_b,F_{cc},F_{cb},F_{bb}$ the partial derivatives of order $1$ and $2$ of $F$
w.r.t.\ the variables $c$ and $b$ (and similarly for $\tilde{F}$).

We need to show that $(c(r,a),b(r,a))$ is infinitely differentiable w.r.t.\ $(r,a)$. To do so, we use
the implicit function theorem. Define
\be{subse}
\cR=\big\{(r,a,c,b)\colon\,(r,a)\in\cL_{\psi^{\hat{\kappa}}},\,
(c,b)\in{\rm int}[\DOM(a)]\big\}
\ee
and
\be{eq18}
\Upsilon_1\colon\,(r,a,c,b)\in \cR\mapsto
(F_{c}+\tilde{F}_{c},F_{b}+\tilde{F}_{b}).
\ee
Let $J_1$ be the Jacobian determinant of $\Upsilon_1$ as a function of $(c,b)$. Applying the implicit
function theorem to $\Upsilon_1$ requires checking three properties:
\begin{itemize}
\item[(i)]
$\Upsilon_1$ is infinitely differentiable on $\cR$.
\item[(ii)]
For all $(r,a) \in \cL_{\psi^{\kappa}}$, the pair $(c(r,a),b(r,a))$ is the only pair in ${\rm int}[\DOM(a)]$
satisfying $\Upsilon_1=0$.
\item[(iii)]
For all $(r,a) \in \cL_{\psi^{\kappa}}$, $J_1 \neq 0$ at $(c(r,a),b(r,a))$.
\end{itemize}
Lemma \ref{smoothkappa} implies that $F$ and $\tilde{F}$ are strictly concave on $\DOM(a)$ and infinitely
differentiable on ${\rm int}[\DOM(a)]$, which is sufficient to prove (i) and (ii). It remains to compute
the Jacobian determinant $J_1$ and prove that it is non-null. This computation is written out
in \cite{dHP07b}, Section 5.4.1, and shows that $J_1$ is non-null when $\tilde{F}_{cc}\tilde{F}_{bb}
-\tilde{F}_{cb}^2>0$. This last inequality is checked in \cite{dHP07b}, Lemma 5.4.2.
\epr

The next step requires an assumption on the set $\cR(p)$. Recall that $\cR_{\alpha,\beta,p}^f$, which is
defined in (\ref{maxdefs}), is the subset of $\cR(p)$ containing the maximizers $(\rho_{kl})$ of the
variational formula in (\ref{fevar}). Consider the triple $(\rho^*(p),\rho^*_{BA}(p),\rho^*_{BB}(p))$,
where $\rho^*(p)$ is defined in (\ref{rho*def}), and
\be{defrho}
\begin{aligned}
\rho^*_{BA}(p) &= \max\{\rho_{BA}\colon\, (\rho_{kl})\in\cR(p) \mbox{ and }
\rho_{AA}+\rho_{AB}=\rho^*(p)\},\\
\rho^*_{BB}(p) &= 1-\rho^*(p)-\rho_{BB}^*(p).
\end{aligned}
\ee

\begin{assumption}
\label{assum}
For all $(\alpha,\beta)\in\cD_2$, $(\rho^*(p),\rho^*_{BA}(p),\rho^*_{BB}(p))
\in\cR_{\alpha,\beta,p}^f$.
\end{assumption}

\noindent
This assumption is reasonable, because in $\cD_2$ (recall Fig.\ \ref{fig-subcrits}) we expect that
the copolymer first tries to maximize the fraction of time it spends crossing $A$-blocks, and then
tries to maximize the fraction of time it spends crossing $B$-blocks that have an $A$-block as
neighbor.

\subsubsection{Smoothness of $f$ on $\cL$}
\label{S5.4}

By Proposition \ref{p:unimaw2}, we know that, for all $r\in(\alpha^*(p),\infty)$, the maximizers
$x(r),y(r),z(r)$ of the variational formula in (\ref{fevarredD2}) are unique. By (\ref{fevar})
and Assumption \ref{assum}, we have that
\be{fevarredD21}
\begin{aligned}
f_{\cD_2}(r) &= V((\rho_{kl}^*(p)),(x(r),y(r),z(r))) = V(\rho^*,x(r),y(r),z(r))\\
&= \frac{\rho^*\, x(r)\, \psi_{AA}(x(r))
+  \rho^*_{BA}\, y(r)\,\psi^{\hat{\kappa}}_{BA}(y(r))
+ \rho^*_{BB}\, z(r)\,\psi_{BB}(z(r))}{\rho^*\, x(r)\,
+ \rho_{BA}^*\, y(r)\,+ \rho_{BB}^*\, z(r)},
\end{aligned}
\ee
where we suppress the $p$-dependence and simplify the notation.

Since $r\in (\alpha^*,\infty)$, Propositions \ref{recobal} and \ref{d2char} imply that $(r,y(r))
\in\cL_{\psi^{\hat{\kappa}}}$. Hence, by Proposition \ref{psiklsmooth}, Corollary \ref{concpsi1}
and the variational formula in (\ref{fevarredD21}), it suffices to prove that $r\mapsto(x(r),y(r),z(r))$
is infinitely differentiable on $(\alpha^*,\infty)$ to conclude that $r\mapsto f_{\cD_2}(r)$ is
infinitely differentiable on $(\alpha^*,\infty)$. To this end, we again use the implicit function
theorem, and define
\be{subse2}
\mathcal{N}=\big\{(r,x,y,z)\colon\, x>2, z>2, (r,y)\in \cL_{\psi^{\hat{\kappa}}} \big\}
\ee
and
\be{eq20**}
\Upsilon_2\colon\,(r,x,y,z)\in\cN \mapsto
\bigg(\frac{\partial V}{\partial x}(\rho^*,x,y,z),
\frac{\partial V}{\partial y}(\rho^*,x,y,z),\frac{\partial V}{\partial z}(\rho^*,x,y,z)\bigg).
\ee
Let $J_2$ be the Jacobian determinant of $\Upsilon_2$ as a function of $(x,y,z)$. To apply the
implicit function theorem, we must check three properties:
\begin{itemize}
\item[(i)]
$\Upsilon_2$ is infinitely differentiable on $\cN$.
\item[(ii)]
For all $r\in(\alpha^*,\infty)$, the triple $(x(r),y(r),z(r))$ is the only triple in
$[2,\infty)^3$ satisfying $(r,x(r),y(r),z(r))\in \cN$ and $\Upsilon_2(r,x(r),y(r),z(r))=0$.
\item[(iii)]
For all $r>\alpha^*$, $J_2 \neq 0$ at $(r,x(r),y(r),z(r))$.
\end{itemize}
Proposition \ref{psiklsmooth} and Corollary \ref{concpsi1}(i) imply that (i) is satisfied. By
Corollary \ref{concpsi1}(ii--iii), we know that $a\mapsto a\psi_{BA}^{\hat{\kappa}}(a)$ and
$a\mapsto a\psi_{kl}(a)$ with $kl\in\{AA,BB\}$ are strictly concave on $[2,\infty)$. Therefore,
Proposition \ref{p:unimaw2} implies that (ii) is satisfied as well. Thus, it remains to prove (iii).

For ease of notation, abbreviate $\psi_{AA}(x)=x\psi_{AA}(r;x)$, $\psi_{BA}(y)=y\psi_{BA}^{\hat{\kappa}}(r;y)$
and $\psi_{BB}(z)=z\psi_{BB}(r;z)$. Note that
\be{secderiv}
\frac{\partial^2 V}{\partial x \partial y}(r,x(r),y(r),z(r))
=\frac{\partial^2 V}{\partial x \partial z}(r,x(r),y(r),z(r))
=\frac{\partial^2 V}{\partial y \partial z}(r,x(r),y(r),z(r))=0,
\ee
which is obtained by differentiating (\ref{fevarredD21}) and using the equality in (\ref{max1}),
i.e.,
\be{secderiv2}
\frac{\partial \psi_{AA}(x)}{\partial x}(x(r))=\frac{\partial \psi_{BA}(y)}{\partial y}(y(r))
=\frac{\partial \psi_{BB}(z)}{\partial z}(z(r))=V(\rho^*,x(r),y(r),z(r)).
\ee
With the help of (\ref{secderiv}), we can assert that, at $(r,x(r),y(r),z(r))$,
\be{21}
J_2=\frac{\partial^2 V}{\partial x^2}\,
\frac{\partial^2 V}{\partial y^2}\,\frac{\partial^2 V}{\partial z^2}
= C\, \frac{\partial^2 \psi_{AA}(x)}{\partial x^2}\,
\frac{\partial^2 \psi_{BA}(y)}{\partial y^2}\,\frac{\partial^2 \psi_{BB}(z)}{\partial z^2},
\ee
with $C$ a strictly positive constant. Abbreviate $x\kappa(x,1)=\kappa(x)$, and denote by
$\kappa''(x)$ its second derivative. Then (\ref{psiAAandBB}) implies that $(\partial^2/\partial x^2)
(\psi_{AA}(x))=\kappa''(x)$ and $(\partial^2/\partial z^2)(\psi_{BB}(z))=\kappa''(z)$. Next, recall
the formula for $\kappa$ stated in \cite{dHW06}, Lemma 2.1.1:
\be{ka}
\kappa(a) = \log 2 + \tfrac12 \left[a \log a - (a-2) \log (a-2)\right],
\qquad a\geq 2.
\ee
Differentiate (\ref{ka}) twice to obtain that $\kappa''$ is strictly negative on $(2,\infty)$.
Therefore it suffices to prove that $(\partial^2/\partial y^2) \psi_{BA}(y(r))<0$ to conclude
that $J_2\neq 0$ at $(r,x(r),y(r),z(r))$, which will complete the proof of Theorem \ref{D2infdiff}.

In \cite{dHP07b}, Lemma 5.5.2, it is shown that the second derivative of $x\psi_{AB}(x)$ w.r.t.\ $x$
is strictly negative at $x^*$, where $(x^*,y^*)$ is the maximizer of the variational formula
in \cite{dHP07b}, Equation (5.5.8), that gives the free energy in the localized phase in the
supercitical regime. It turns out that this proof readily extends to our setting, and for this
reason we do not repeat it here.

\subsection{Proof of Theorem \ref{th:orrder}}
\label{S4.3}

\bpr
Recall that $\alpha^*(p)=\alpha^*$ and set $\rho^*(p)=\rho^*$. Let $x_\delta,y_\delta$ be
the unique maximizers of the variational formula in (\ref{fevarred}) at $\alpha^*+\delta$,
i.e.,
\be{fevarred7}
f_{\cD_1}(\alpha^{*}+\delta;p) =
\frac{\rho^{*}\,x_\delta\,\kappa(x_\delta,1)+(1-\rho^{*})\,y_\delta\,
[\kappa(y_\delta,1)-\frac12(\alpha^{*}+\delta)]}
{\rho^{*}\,x_\delta + (1-\rho^{*})\,y_\delta}.
\ee
Put
\be{Tdeltadef}
T_\delta = f_{\cD_2}(\alpha^*+\delta)-f_{\cD_1}(\alpha^*)
-f'_{\cD_1}(\alpha^*)\,\delta -\tfrac12\,f''_{\cD_1}(\alpha^*)\,\delta^2
\ee
and $V_{\delta}=\rho^* x_\delta+(1-\rho^*)y_\delta$. By picking $x=x_\delta$, $y=y_\delta$
and $z=y_\delta$ in (\ref{fevarredD2}), we obtain that, for every $(b,c)\in \DOM(y_\delta)$,
\be{eq:borninf}
\begin{aligned}
f_{\cD_2}(\alpha^*+\delta)\geq \frac{1}{V_\delta}
\bigg(\rho^{*} x_\delta\, \kappa(x_\delta,1)&+\rho_{BA}\bigg[c
\hat{\kappa}(\tfrac{c}{b})+(y_\delta-c)\Big[\kappa\Big(y_\delta-c,1-b\Big)
-\tfrac{\alpha^*+\delta}{2}\Big]\bigg]\\
&+(1-\rho^*-\rho_{BA})\,y_\delta\, \Big[\kappa(y_\delta,1)-\tfrac{\alpha^*+\delta}{2}\Big]\bigg).
\end{aligned}
\ee
Hence, using a first-order Taylor expansion of $(a,b)\mapsto\kappa(a,b)$ at $(y_\delta,1)$, and
noting that $b\leq c$ and $V_\delta\geq 2$ for all $\delta\geq 0$, we obtain
\be{eq:borninf2}
f_{\cD_2}(\alpha^*+\delta)
\geq f_{\cD_1}(\alpha^{*}+\delta)
+\frac{\rho_{BA}}{(c/b)\,V_\delta}\, R_{\tfrac{c}{b},\delta}\,c + L(\delta)\, c^2,
\ee
where $\delta\mapsto L(\delta)$ is bounded in the neighborhood of $0$ and
\be{Rmudel}
R_{\mu,\delta} = \Big[\mu
\Big(\hat{\kappa}(\mu)-\kappa(y_\delta,1)-y_\delta\,\partial_1\kappa(y_\delta,1)
+\tfrac{\alpha^*+\delta}{2}\Big)-y_\delta \partial_2\kappa(y_\delta,1)\Big].
\ee
The strict concavity of $\mu\mapsto \mu \hat{\kappa}(\mu)$ implies that, for every $\delta>0$,
$\mu \mapsto R_{\mu,\delta}$ attains its maximum at a unique point $\mu_\delta$. Thus, we may
pick $b=c/\mu_0$ in (\ref{eq:borninf2}) and obtain
\be{eq:order}
T_\delta\geq \Big\{f_{\cD_1}(\alpha^{*}+\delta)-f_{\cD_1}(\alpha^*)-f'_{\cD_1}(\alpha^*)\,\delta
-\frac{1}{2}\,f''_{\cD_1}(\alpha^*)\,\delta^2\Big\}
+\frac{\rho_{BA}}{\mu_0\,V_\delta}\, R_{\mu_0,\delta}\,c
+ L(\delta)\, c^2.
\ee
Since $(\alpha,\beta)\mapsto f_{\cD_1}(\alpha-\beta)$ is analytic on $\CONE$, and since $\delta
\mapsto V_\delta$ is continuous (recall that, by Proposition \ref{p:unimaw1}, $\delta\mapsto
(x_\delta,y_\delta)$ is continuous), we can write, for $\delta$ small enough,
\be{eq:order2alt}
T_\delta\geq \frac{\rho_{BA}}{2\mu_0\,V_0}\, R_{\mu_0,\delta}\,c
+L(\delta)\, c^2+G(\delta)\,\delta^3
\ee
where $\delta\mapsto G(\delta)$ is bounded in the neighborhood of $0$. Next, note that
Proposition~\ref{philbextext} implies that $R_{\mu_0,0}=0$. Moreover, as shown in \cite{dHW06},
Proposition 2.5.1, $\delta\mapsto (x_\delta,y_\delta)$ is infinitely differentiable, so that
$\tfrac{\partial R}{\partial \delta}(\mu_0,0)$ exists. If the latter is $>0$, then we pick
$c=x \delta$ in \eqref{eq:order2} and, by choosing $x>0$ small enough, we obtain that there
exists a $t>0$ such that $T_\delta\geq t \delta^2$ and the proof is complete.

Thus, it remains to prove that $\tfrac{\partial R}{\partial \delta}(\mu_0,0)>0$. To that aim,
we let $(x'_0, y'_0)$ be the derivative of $\delta\mapsto (x_\delta,y_\delta)$ at
$\delta=0$ and recall the following expressions from \cite{dHW06}:
\be{kappaform}
\begin{aligned}
\kappa(a,1)
&= \tfrac{1}{a}\Big[\log 2+\tfrac{1}{2}\,\big[a\log a-(a-2)\log(a-2)\big]\Big],\\
\frac{\partial\kappa}{\partial 1}(a,1)
&=-\tfrac{\log 2}{a^2}-\tfrac{1}{a^2}\log(a-2),\\
\frac{\partial\kappa}{\partial 2}(a,1)
&=\tfrac{1}{2a} \log\Big(\tfrac{4 (a-2)(a-1)^2}{a}\big).
\end{aligned}
\ee
These give
\be{dude1}
\frac{\partial R}{\partial \delta}(\mu_0,0)
=\frac{1}{2}\bigg\{\mu_0+y'_0\Big[\frac{2 \mu_0}{(y_0-2)y_0}
-\frac{1}{y_0-2}-\frac{2}{y_0-1}+\frac{1}{y_0}\Big]\bigg\}.
\ee
Since $-\frac{1}{y_0-2}-\frac{2}{y_0-1}+\frac{1}{y_0}<0$, and since it was proven
in \cite{dHW06}, Proposition 2.5.1, that $y'_0<0$, we obtain via (\ref{dude1}) that
\be{dude2}
\frac{\partial R}{\partial \delta}(\mu_0,0)
\geq \frac{\mu_0}{2}\Big(1+y'_0\frac{2}{(y_0-2)y_0}\Big).
\ee
It was also proven in \cite{dHW06}, Proposition 2.5.1, that
\be{y_0}
y'_0=-\frac{y_0(y_0-2)}{2}+\frac{2 x'_0 y_0 (y_0-2)}{x_0 (x_0-2)},
\ee
which implies that $y'_0>-\tfrac{y_0(y_0-2)}{2}$ because $x'_0>0$. Thus, recalling
(\ref{dude2}), we indeed have that $\frac{\partial R}{\partial \delta}(\mu_0,0)>0$.
\epr

\subsection{Proof of Theorem \ref{th:orrder2}}
\label{S4.4}

\subsubsection{Lower bound}
\label{S4.4.1}

\bpr
Pick $r\in[0,\alpha^*)$ and $\delta>0$. Denote by $\alpha_r$ and $\beta_r$ the quantities
$r+\beta_c^1(r)$ and $\beta_c^1(r)$. Let $x_r,y_r$ be the maximizers of (\ref{fevarred})
at $\alpha_r-\beta_r$ (keep in mind that $\alpha_r-\beta_r=r$), i.e.,
\be{fevarred7alt}
f_{\cD_1}(\alpha_r-\beta_r;p)
= \frac{\rho^{*}\, x_r\, \kappa(x_r,1)+(1-\rho^{*})\, y_r\,
[\kappa(y_r,1)-\frac{\alpha_r-\beta_r}{2}]}{\rho^{*} \,x_r + (1-\rho^{*})\, y_r}.
\ee
Put
\be{Tdeltadefalt}
T_\delta=f_{\cL_1}(\alpha_r+\delta,\beta_r+\delta;p)-f_{\cD_1}(\alpha_r-\beta_r;p)
\ee
and $V_{r}=\rho^* x_r+(1-\rho^*)y_r$. By picking $x=x_r$, $y=y_r$ and $z=y_r$ in (\ref{fevarredL1})
at $(\alpha_r+\delta,\beta_r+\delta)$, we obtain, for every $(b,c)\in \DOM(y_r)$,
\be{eq:borinf}
\begin{aligned}
&f_{\cL_1}(\alpha_r+\delta,\beta_r+\delta;p)\\
&\qquad\geq \frac{1}{V_r}\bigg(\rho^{*} x_r\, \kappa(x_r,1)
+\rho_{BA}\bigg[c\phi^{\cI}(\alpha_r+\delta,\beta_r+\delta;\tfrac{c}{b})
+(y_r-c)\Big[\kappa\Big(y_r-c,1-b\Big)
-\tfrac{\alpha_r-\beta_r}{2}\Big]\bigg]\\
&\qquad\qquad +(1-\rho^*-\rho_{BA})\,y_r\,
\Big[\kappa(y_r,1)-\tfrac{\alpha_r-\beta_r}{2}\Big]\bigg).
\end{aligned}
\ee
Therefore, using \eqref{fevarred7alt} and \eqref{eq:borinf}, we obtain
\be{eq:order21}
T_\delta\geq \frac{\rho_{BA}}{V_r}
\bigg(c\Big[\phi^{\cI}(\alpha_r+\delta,\beta_r+\delta;\tfrac{c}{b})
+ \tfrac{\alpha_r-\beta_r}{2}-\kappa(y_r-c,1-b)\Big]
+ y_r\Big[\kappa(y_r-c,1-b)-\kappa(y_r,1)\Big]\bigg).
\ee
By a Taylor expansion of $(a,b)\mapsto \kappa(a,b)$ at $(y_r,1)$, noting that $b\leq c$,
we can rewrite
(\ref{eq:order21}) as
\be{eq:order21alt}
T_\delta\geq \frac{\rho_{BA}}{V_r}
\bigg(c\Big[\phi^{\cI}(\alpha_r+\delta,\beta_r+\delta;\tfrac{c}{b})
+ \tfrac{\alpha_r-\beta_r}{2}-\kappa(y_r,1)-y_r \partial_1\kappa(y_r,1)
-\tfrac{b\, y_r}{c}\partial_2\kappa(y_r,1)\Big]+\xi(c) c^2\bigg),
\ee
where $x\mapsto \xi(x)$ is bounded in the neighborhood of $0$. As explained in the proof of
Theorem \ref{phtrcurve}(i), for $r<\alpha^*$,
\be{maxim2}
\sup_{\mu\geq 1}\mu \Big(\hat{\kappa}(\mu)
+ \tfrac{\alpha_r-\beta_r}{2}-\kappa(y_r,1)-y_r \partial_1\kappa(y_r,1)-\tfrac{y_r}{\mu}
\partial_2\kappa(y_r,1)\Big)<0.
\ee
Set
\be{sr}
s_r=\tfrac{\alpha_r-\beta_r}{2}-\kappa(y_r,1)-y_r \partial_1\kappa(y_r,1).
\ee
Then \eqref{maxim2} and \cite{dHW06}, Lemma 2.1.2(iii), which asserts that $\mu \hat{\kappa}(\mu)
\sim \log\mu$ as $\mu \to \infty$, are sufficient to conclude that $s_r<0$. Next, we note that,
by the definition of $\beta_r$ and by Corollary \ref{philbextext}, for all $\delta>0$ there
exists a $\mu_\delta>1$ such that
\be{bilout}
\mu_\delta \phi^{\cI}(\alpha_r+\delta,\beta_r+\delta;\mu_\delta)
+ \mu_\delta s_r -y_r \partial_2\kappa(y_r,1)>0.
\ee
Because of Lemma \ref{mulim}(i), which tells us that $\phi^{\cI}(\alpha_r+\delta,\beta_r+\delta;\mu)$
tends to $0$ as $\mu\to\infty$ uniformly in $\delta\in[0,1]$, we know that $\mu_\delta$ is necessarily
bounded uniformly in $\delta$. For this reason, and since $(\alpha_r,\beta_r)\in \cD_1$ and
$(\mu,\alpha,\beta)\mapsto \phi^{\cI}(\alpha,\beta;\mu)$ is continuous, Corollary \ref{philbextext}
allows us to assert that there exists a $\mu_r>1$ such that
\be{maxim}
\mu_r \Big(\phi^{\cI}(\mu_r,\alpha_r,\beta_r)
+s_r-\tfrac{y_r}{\mu_r} \partial_2\kappa(y_r,1)\Big) = 0.
\ee
Hence, using (\ref{maxim2}), we obtain $\phi^{\cI}(\alpha_r,\beta_r,\mu_r)>\hat{\kappa}(\mu_r)$.
Moreover, $x\mapsto \phi^{\cI}(\alpha_r-\beta_r+x,x,\mu_r)$ is convex and $\phi^{\cI}(\alpha_r-\beta_
r+x,x,\mu_r)=\hat{\kappa}(\mu_r)$ for $x\leq 0$. Therefore, we can assert that $\partial_x(\phi^{\cI}
[\alpha_r-\beta_r+x,x,\mu)](x=\beta_r)=m>0$ and, consequently, $\phi^{\cI}(\mu_r,\alpha_r+\delta,
\beta_r+\delta)\geq \phi^{\cI}(\mu_r,\alpha_r,\beta_r)+m \delta$. Now (\ref{eq:order21alt}) becomes
\be{eq:order22}
T_\delta\geq \frac{\rho_{BA}}{V_r} \big(c m \delta+\xi(c) c^2\big),
\ee
and by picking $c=\sigma \delta$ with $\sigma$ small enough we get the claim.
\epr

\subsubsection{Upper bound}
\label{S4.4.2}

\bpr
For this proof only we assume the strict concavity of $\mu \mapsto \mu \phi^{\cI}(\alpha,\beta;\mu)$
for all $(\alpha,\beta)\in \CONE$. As mentioned in Remark \ref{strconcav}(ii), the latter implies
the strict concavity of $a\mapsto a\psi_{BA}(\alpha,\beta;a)$. We keep the notation of
Section~\ref{S4.4.1}, i.e., we let $(x_\delta,y_\delta,z_\delta)$ be the unique maximizers of the
variational formula in (\ref{fevarredL1}) at $(\alpha_r+\delta,\beta_r+\delta)$:
\be{rep}
f_{\cL_1}(\alpha_r+\delta,\beta_r+\delta;p)
= \frac{\rho_A\, x_\delta\, \psi_{AA}(x_\delta)
+  \rho_{BA}\, y_\delta\,\psi_{BA}(y_\delta) + \rho_{BB}\, z_\delta\,\psi_{BB}(z_\delta)}
{\rho_A\, z_\delta\, + \rho_{BA}\, y_\delta\,+ \rho_{BB}\, z_\delta}.
\ee
We let $(b_\delta,c_\delta)$ be the maximizers of the variational formula in \eqref{psiinflink}
at $(\alpha_r+\delta,\beta_r+\delta,z_\delta)$. Recall that, in $\cD_1$, $f(\alpha,\beta;p)$
is equal to $f_{\cD_1}(\alpha,\beta;p)$, but is also equal to $f_{\cL_1}(\alpha,\beta;p)$.
Therefore, by picking $x_\delta$, $y_\delta$ and $z_\delta$ in (\ref{fevarredL1}) at
$(\alpha_r,\beta_r)$, and $(b_{\delta},c_\delta)$ in (\ref{psiinflink}) at $(\alpha_r,
\beta_r,y_\delta)$, we obtain the upper bound
\be{maj}
\begin{aligned}
T_\delta &= f_{\cL_1}(\alpha_r+\delta,\beta_r+\delta,p)
- f_{\cL_1}(\alpha_r,\beta_r,p)\\
&\leq c_{\delta} \big[\phi^{\cI}(\alpha_r+\delta,\beta_r+\delta,
\tfrac{c_\delta}{b_\delta})-\phi^{\cI}(\alpha_r,\beta_r,\tfrac{c_\delta}{b_\delta})\big].
\end{aligned}
\ee
As stated in the proof of Lemma \ref{smo}, $\phi^{\cI}(\alpha,\beta;\mu)$ is $K$-Lipshitz in
$(\alpha,\beta)$ uniformly in $\mu$. We therefore deduce from \eqref{maj} that $T_\delta \leq
K c_{\delta} \delta$, and the proof will be complete once we show that $\lim_{\delta\downarrow 0}
c_{\delta}=0$.

Since $(\alpha_r,\beta_r)\in \cD_1$, we have $f_{\cD_1}=f_{\cL_1}$ at $(\alpha_r,\beta_r)$ and
since $a\mapsto a \psi_{BA}(\alpha_r,\beta_r;a)$ is strictly concave, the maximizers of
\eqref{fevarredL1} are unique. This yields $y_0=z_0$. Moreover, by applying Proposition
\ref{p:unimaw1}, we obtain that $(x_\delta,z_\delta)\mapsto (x_0,z_0)$ as $\delta\to 0$. Therefore
we need to show that $y_\delta \mapsto y_0$ as $\delta\to 0$. For this, we recall \eqref{max1},
which allows us to assert that, for $\delta\geq 0$,
\be{ret}
\frac{\partial [a\psi_{BA}(\alpha_r+\delta,\beta_r+\delta;a)]}{\partial a}(y_\delta)
=f(\alpha_r+\delta,\beta_r+\delta).
\ee
Since $f(\alpha_r,\beta_r)>0$, by applying Lemma \ref{l:c'}(ii) we can assert that there exists
an $a_0>0$ such that, for all $\delta\in[0,1]$ and all $y\geq a_0$, the l.h.s.\ of \eqref{ret}
is smaller than or equal to $f(\alpha_r,\beta_r)/2$, whereas for $\delta$ small enough the
continuity of $f$ implies that the r.h.s.\ of \eqref{ret} is strictly larger than
$f(\alpha_r,\beta_r)/2$. Therefore $y_\delta\leq a_0$ when $\delta$ is small. Next, the continuity
of $(\alpha,\beta)\mapsto f(\alpha,\beta)$ and $(\alpha,\beta,a)\mapsto \psi_{kl}(\alpha,\beta;a)$,
together with the convergence of $(x_\delta,z_\delta)$ to $(x_0,y_0)$ as $\delta\to 0$ allows us,
after letting $\delta\to 0$ in \eqref{rep} and using again the uniqueness of the maximizers in
\eqref{fevarredL1}, to conclude that $y_\delta\to y_0$ as $\delta \to 0$.

Now, put $\mu_{\delta}=\tfrac{c_\delta}{b_\delta}$. Then, by the definition of $(b_\delta,c_\delta)$,
we have
\be{rel}
y_\delta\, \psi_{BA}(\alpha_r+\delta,\beta_r+\delta;y_\delta)
=c_\delta \phi^{\cI}(\alpha_r+\delta,\beta_r+\delta;\mu_\delta)
+(y_\delta-c_\delta) \Big[\kappa(y_\delta-c_\delta,1-\tfrac{c_\delta}{\mu_\delta})
-\tfrac{\alpha_r-\beta_r}{2}\Big].
\ee
We already know that $c_\delta$ is bounded in $\delta$, since $c_\delta\leq y_\delta$ and $y_\delta$
converges. We want to show that $\mu_\delta$ is bounded in $\delta$ as well. Note that, by the
concavity of $a\mapsto a\kappa(a,1)$, the r.h.s.\ of \eqref{rel} is concave as a function of $c$.
Moreover, for all $\delta>0$ we have $c_\delta>0$, because $(\alpha_r+\delta,\beta_r+\delta)\in\cL_1$.
This implies that the derivative of the r.h.s.\ of \eqref{rel} w.r.t.\ $c$ at $(c=0, \mu_\delta)$
is strictly positive when $\delta>0$, i.e.,
\be{rev}
\phi^{\cI}(\alpha_r+\delta,\beta_r+\delta;\mu_\delta)
-s_{r,\delta}-\tfrac{y_\delta}{\mu_\delta} \partial \kappa(y_\delta,1)>0
\ee
with
\be{N13}
s_{r,\delta}=\tfrac{\alpha_r-\beta_r}{2}-\kappa(y_{\delta},1)-y_{\delta} \partial_1\kappa(y_{\delta},1).
\ee
Note that $s_{r,0}=s_r$, with $s_r$ defined in \eqref{sr}. Since $y_\delta$ converges to $y_0$, and
since we have proved in Section \ref{S4.4.1} that $s_{r,0}<0$, it follows that $s_{r,\delta}\leq
s_{r,0}/2<0$ for $\delta$ small enough. Moreover, by Lemma \ref{mulim}(i), $\phi^{\cI}(\alpha_r+\delta,
\beta_r+\delta;\mu)$ tends to $0$ as $\mu\to \infty$ uniformly in $\delta\in [0,1]$, which, with the
help of \eqref{rev}, is sufficient to assert that $\mu_\delta$ is bounded from above for $\delta$
small.

At this stage, it remains to prove that the only possible limit for $c_\delta$ is $0$. Assuming
that $(c_\delta,\mu_\delta)\mapsto (c_\infty,\mu_\infty)$, we obtain, when $\delta \to 0$ in
\eqref{rel},
\be{rel1}
y_0\, \psi_{BA}(\alpha_r,\beta_r;y_0)=c_\infty \phi^{\cI}(\alpha_r,\beta_r;\mu_\infty)+
(y_0-c_\infty) \Big[\kappa(y_0-c_\infty,1-\tfrac{c_\infty}{\mu_\infty})
-\tfrac{\alpha_r-\beta_r}{2}\Big].
\ee
The fact that $(\alpha_r,\beta_r)\in\cD_1$ implies, by Corollary \ref{philbextext}, that the
derivative of the r.h.s.\ of \eqref{rel1} w.r.t.\ $c$ at $(c=0,\mu_\infty)$ is non-positive.
Therefore the concavity in $c$ of the r.h.s.\ of \eqref{rel1} is sufficient to assert that
$c_\infty=0$.
\epr

\subsection{Proof of Theorem \ref{th:orrder3}}
\label{S4.5}

\bpr
Recall Theorem \ref{phtrcurve2}(iv), and the constant $r_2>0$ such that $\cD_2$ and $\cL_1$ touch
each other along the curve $r\in[\alpha^{*},\alpha^*+r_2)\mapsto (r+\beta_c^{2}(r),\beta_c^2(r))$.
Pick $r\in[\alpha^*,\alpha^*+r_2)$ and $\delta>0$. We abbreviate $\alpha_r$ and $\beta_r$ for the
quantities $r+\beta_c^2(r)$ and $\beta_c^2(r)$. Let $x_r,y_r,z_r$ be the unique maximizers of the
variational formula \eqref{fevarredD2} at $\alpha_r-\beta_r$, i.e.,
\be{fevarred8alt}
f_{\cD_2}(\alpha_r-\beta_r;p)
= \frac{\rho_A\, x_r\, \psi_{AA}(\alpha_r-\beta_r;x_r)+\rho_{BA}\, y_r\,
\psi_{BA}^{\hat{\kappa}}(\alpha_r-\beta_r;y_r)+\rho_{BB} z_r
\psi_{BB}(\alpha_r-\beta_r;z_r)}{\rho_{A} x_r+\rho_{BA} y_r +
\rho_{BB} z_r}.
\ee
Put
\be{Tdeltadefalt2}
T_\delta=f_{\cL_1}(\alpha_r+\delta,\beta_r+\delta;p)-f_{\cD_2}(\alpha_r-\beta_r;p)
\ee
and $V_{r}=\rho_{A} x_r+\rho_{BA} y_r +
\rho_{BB} z_r$. By picking $x=x_r$, $y=y_r$ and $z=y_r$ in \eqref{fevarredL1},
we obtain,
\begin{align}\label{impL1}
\nonumber f_{\cL_1}(\alpha_r+\delta,&\beta_r+\delta;p)\geq\\
&\frac{\rho_A\, x_r\, \psi_{AA}(\alpha_r-\beta_r,x_r)+\rho_{BA}\, y_r\,
\psi_{BA}(\alpha_r+\delta,\beta_r+\delta;y_r)+\rho_{BB} z_r \psi_{BB}(\alpha_r-\beta_r;z_r)}{V(r)}.
\end{align}
Therefore, using (\ref{fevarred8alt}--\ref{impL1}) we obtain
\be{diffi}
T_\delta\geq \frac{\rho_{BA} y_r}{V_r}
\Big(\psi_{BA}(\alpha_r+\delta,\beta_r+\delta; y_r)
-\psi_{BA}^{\hat{\kappa}}(\alpha_r-\beta_r; y_r)\Big).
\ee
Let $(c_r,b_r)$ be the unique maximizer of \eqref{psiinflink2} at $(\alpha_r-\beta_r;y_r)$. By
picking $(c,b)=(c_r,b_r)$ in \eqref{psiinflink} at $(\alpha_r+\delta,\beta_r+\delta;y_r)$,
we can bound $T_\delta$ from below as
\be{eq:order23alt}
T_\delta\geq R_r
\bigg(\Big[\phi^{\cI}(\alpha_r+\delta,\beta_r+\delta;\tfrac{c_r}{b_r})
-\hat{\kappa}(\tfrac{c_r}{b_r})\Big]\bigg),
\ee
where $R_r=\rho_{BA}c_r/V_r$. Since $(\alpha_r,\beta_r)\in\cD_2$ and $(\alpha_r+\delta,\beta_r+\delta)
\in \cL_1$, it follows from Proposition \ref{d2char} that $\psi_{BA}(\alpha_r,\beta_r;y_r)
=\psi_{BA}^{\hat{\kappa}}(\alpha_r,\beta_r;y_r)$ and $\psi_{BA}(\alpha_r+\delta,\beta_r+\delta;y_r)
>\psi_{BA}^{\hat{\kappa}}(\alpha_r,\beta_r;y_r)$. Therefore, by Lemma \ref{linint}, we obtain that
$\phi^{\cI}(\alpha_r,\beta_r;\tfrac{c_r}{b_r})=\hat{\kappa}(\tfrac{c_r}{b_r})$, whereas
$\phi^{\cI}(\alpha_r+\delta,\beta_r+\delta;\tfrac{c_r}{b_r})>\hat{\kappa}(\tfrac{c_r}{b_r})$,
which means that the phase transition of $\phi^{\cI}$ along $\{(s+r,s)\colon\; s\geq -\tfrac{r}{2}\}$
effectively occurs at $s=\beta_c(r)$. Using \eqref{th:orrder3}, we complete the proof.
\epr

%%%%%%%%%%%%%%%%%%%%%%%%%%%%% REFERENCES %%%%%%%%%%%%%%%%%%%%%%%%%%%%%%%%%%%

\end{document}